\documentclass[final]{siamart251216}
\ifpdf
\hypersetup{
  pdftitle={Stochastic Galerkin Method and Hierarchical Preconditioning for PDE-constrained Optimization},
  pdfauthor={Z. Li, A. Onwunta, and B. Soused\'{\i}k}
}
\fi
\newtheorem{remark}{Remark}
\usepackage{microtype}
\usepackage{algorithm}
\usepackage{algpseudocode}
\usepackage{amsfonts}
\usepackage{amsmath}
\usepackage{amssymb}
\usepackage{enumerate}
\usepackage{eurosym}
\usepackage{mathrsfs}
\usepackage{mathtools}
\usepackage{diagbox}
\usepackage{graphicx}
\usepackage{multirow}
\usepackage{url}
\usepackage{verbatim}
\usepackage{xcolor}

\newcommand{\zhl}[1]{{\color{black}#1}}
\newcommand{\zh}[1]{{\color{black}#1}}

\usepackage[pdftex]{thumbpdf}
\providecommand{\U}[1]{\protect\rule{.1in}{.1in}}
\DeclareMathOperator{\vectorize}{vec}
\begin{document}

% alternative title:
% Hierarchical Preconditioning for PDE-based Optimal Control Problems with Uncertain Coefficients
\headers{A Preconditioner for Stochastic Optimal Control}{Z. Li, A. Onwunta, and B. Soused\'{\i}k}

\title{Stochastic Galerkin Method and Hierarchical Preconditioning for PDE-constrained Optimization
    \thanks{Dedicated to the memory of Professor Howard C. Elman, our dear colleague and mentor.
        We also want to thank Lehigh's High Performance Computing systems for providing computational resources.
    }
    %\thanks{This work is based upon work
    %supported by the U.\, S.\, Department of Energy Office of Advanced Scientific
    %Computing Research, Applied Mathematics program under Award Number
    %DE-SC0009301, and by the U.\, S.\, National Science Foundation under grants
    %DMSâ1418754 and DMS1521563.}
}
\author{Zhendong Li\thanks{Department of Industrial and Systems Engineering,
    Lehigh University, Bethlehem, PA 18015 (\texttt{zhl923@lehigh.edu}, \texttt{ako221@lehigh.edu}).}
\and Akwum Onwunta\footnotemark[2]
\and Bed\v{r}ich Soused\'{\i}k\thanks{Department of Mathematics and Statistics,
    University of Maryland, Baltimore County, 1000 Hilltop Circle, Baltimore,
    MD~21250 (\texttt{sousedik@umbc.edu}). }}
\maketitle
\begin{abstract}
    We develop efficient hierarchical preconditioners for optimal control problems governed by partial differential equations with uncertain coefficients. Adopting a discretize-then-optimize framework that integrates finite element discretization, stochastic Galerkin projection, and advanced time-discretization schemes, the approach addresses challenges of scaling large and ill-conditioned linear systems arising in uncertainty quantification. By exploiting sparsity of linear systems in stochastic Galerkin method, we formulate hierarchical preconditioners based on truncated stochastic expansion that strike an effective balance between computational cost and preconditioning quality. Numerical experiments demonstrate that the proposed preconditioners significantly accelerate the convergence of iterative solvers compared to existing methods, providing robust and efficient solvers for both steady-state and time-dependent optimal control problems under uncertainty.
\end{abstract}
% REQUIRED

% REQUIRED

\section{Introduction}
\label{sec:introduction}

Optimal control problems governed by partial differential~equations (PDEs) arise in numerous applications,
including fluid mechanics, %~\cite{antil2012reduced,antil2015optimal}, 
structural optimization, %~\cite{lohner2020detailed}, 
and inverse problems. %~\cite{antil2023deep}. 
These problems have been extensively studied over the past decades.
For a theoretical overview and computational methods related to deterministic problems, we refer readers to, e.g.,~\cite{elman_finite_2014,Troltzsch-2010-OCP}.
In many practical applications, the PDE coefficients are uncertain.
Such uncertainties originate from various sources, including measurement errors, model approximations, and environmental variations,
and they are modeled as random variables or stochastic processes.
Recently, there has been an increased interest in optimal control problems governed by PDEs with random coefficients, see e.g.,~\cite{Martinez-Frutos-2018-OCP}.
These stochastic problems are inherently more complex than their deterministic counterparts, thus necessitating specialized numerical methods.

% Two alternative strategies are used for the optimal control problems: optimize-then-discretize and discretize-then-optimize.
% The optimize-then-discretize approach involves deriving continuous optimality conditions and then discretizing them. 
% %This method is often preferred for theoretical studies, enabling rigorous analysis of optimality conditions and numerical convergence. 
%  Conversely, the discretize-then-optimize approach discretizes the objective and the PDE first and then solves the resulting discrete optimization problem. The discretize-then-optimize is widely used in practice, because it allows the employment of efficient numerical methods 
% such as finite element or finite difference methods, see e.g., ~\cite{HDO2023, BDOS2016,  Benner-2017-SOC,Benner-2016-BDP, elman_finite_2014}. 
% This approach is also used in this study. 
Two strategies are commonly used: optimize-then-discretize and discretize-then-optimize. We adopt the latter, which discretizes the objective and PDE first, enabling the use of efficient numerical methods like finite elements~\cite{HDO2023, BDOS2016, Benner-2017-SOC, Benner-2016-BDP, elman_finite_2014}.

% In the formulation of numerical methods, discretizations in the stochastic space, spatial domain, and potentially time domain are required. 
% For stochastic-space discretization, several methods exist, including Monte Carlo method, 
% stochastic collocation method and stochastic Galerkin method. 
% The stochastic collocation method discretizes random variables using a set of collocation points, 
% solving the resulting deterministic PDEs at these points~\cite{Chen-2014-WRB,Kouri-2013-MSC,Kouri-2013-TRA,Tiesler-2012-SCO}. 
% The stochastic Galerkin method 
% expands the solution in terms of orthogonal polynomials and solves the resulting coupled system of equations ~\cite{Ghanem-1991-SFE, LeMaitre-2010-SMU, Lee-2013-SGM, Lord-2014-ICS,  Rosseel-2012-OCS, Xiu-2010-NMS}. 
% %In this paper, we primarily focus on the stochastic Galerkin method.

% Monte Carlo method is simple and therefore popular. However, typically a large number of samples is required for acceptable accuracy, making it computationally intensive for high-dimensional problems. 
% The stochastic collocation method is in general more efficient, but it may lose efficiency in high-dimensional random spaces, 
% and it becomes challenging to implement in practice. 
% Here, we use the stochastic Galerkin method, which systematically captures uncertainties by employing orthonormal polynomial expansions. It thus 
% preserves optimal convergence properties, and a careful implementation allows to improve computational efficiency for large-scale problems.
For stochastic discretization, we employ the stochastic Galerkin method~\cite{Ghanem-1991-SFE, LeMaitre-2010-SMU, Lee-2013-SGM, Lord-2014-ICS, Rosseel-2012-OCS, Xiu-2010-NMS}, which expands the solution using orthonormal polynomials. Although the Monte Carlo method is simpler and the stochastic collocation method~\cite{Chen-2014-WRB, Kouri-2013-MSC} decouples the problem, stochastic Galerkin method systematically captures uncertainties and preserves optimal convergence properties, offering computational efficiency for large-scale problems when carefully implemented.

% % Specifically, let 
% % ${t_0,t_1,\dots,t_{n_t}}$ denote a partition of the time interval of interest, 
% % ${\phi_i(\mathbf{x})}_{i=1}^{n_h}$ represent a spatial discretization basis
% % and ${\psi_\ell({\xi})}_{\ell=1}^{n_\xi}$ represent an orthonormal polynomial basis in the random space. 
% % Then the state $y(t,\mathbf{x},{\xi})$ can be approximated by
% % \[
% % y(t,\mathbf{x},{\xi}) \approx 
% % \sum_{k=0}^{n_t}\sum_{i=1}^{n_h}\sum_{\ell=1}^{n_\xi} y_{k,i,\ell}\,\Theta_k(t)\,\phi_i(\mathbf{x})\,\psi_\ell({\xi}),
% % \]
% % where $\Theta_k(t)$ is a temporal basis, and $y_{n,i,k}$ are expansion coefficients (i.e., the degrees of freedom of the numerical solution) in the basis $\{\Theta_k(t)\,\phi_i(\mathbf{x})\,\psi_\ell({\xi})\}$, which are determined by a suitable numerical method.

In order to address the PDE-constrained optimization problems, we employ the discretize-then-optimize approach.
In practice, it is common to combine temporal discretization (e.g., backward Euler scheme), stochastic expansions (e.g., generalized polynomial chaos expansions), and spatial discretization (e.g., finite element method).
This leads to large-scale linear systems obtained via finite element (or possibly finite difference) discretizations,
which are then typically solved using Krylov subspace methods, e.g., by the generalized minimal residual method (GMRES).
These linear systems are often ill-conditioned, which causes a slow convergence of the iterative method.
To address this issue, we construct preconditioners that improve the convergence of iterative solvers.
% We note that the development of efficient solvers is a general challenge in optimal control, including for deterministic problems solved via duality-based approaches~\cite{MR4750911}. 

In this paper, we introduce a hierarchical preconditioning framework specifically tailored for stochastic PDE-constrained optimal control problems.
Although the core concepts are inspired by preconditioners for PDE problems~\cite{Bespalov-2021-TPS,Sousedik-2014-THP,Sousedik-2014-HSC,Youngnoi-2025-ICS}, an application to the Karush-Kuhn-Tucker (KKT) systems arising from optimization problems and an extension to time-dependent problems are nontrivial
and constitute a primary contribution of this work.
We provide a systematic derivation of the preconditioner for both steady-state and all-at-once formulation of time-dependent problems.
The method is supported by a spectral analysis,
proving that the proposed preconditioner is spectrally equivalent to the ideal but computationally prohibitive exact Schur complement.
%which theoretically guarantees its robustness and scalability. 
The performance is evaluated using a set of numerical experiments.

The paper is organized as follows.
In Section~\ref{sec:steady-state problem}, we introduce the steady-state stochastic optimal control problem and its discretization into a large-scale KKT system.
Then, we construct the hierarchical Gauss-Seidel preconditioner for PDE-constrained optimal control problems (hGSoc).
In Section~\ref{sec:Time-Dependent Problem}, we extend the framework to the time-dependent case.
We present an all-at-once  discretization and develop a corresponding parallel-in-time preconditioner.
In Section~\ref{sec: spectral analysis}, we provide a spectral analysis of the preconditioners.
In Section~\ref{sec:numerical-experiments}, we demonstrate the efficiency of the methods
by a set of numerical experiments.
Finally, in Section~\ref{sec:conclusion} we conclude and summarize our work.
\section{Steady-state problem}
\subsection{Problem formulation}
\label{sec:steady-state problem}
% Let $\left(\Omega,\mathcal{F},\mathcal{P}\right)$ be a complete probability space, where~$\Omega$ is the sample space, 
% $\mathcal{F}$ is the $\sigma$-algebra generated by $\Omega$, and $\mathcal{P}$ is the probability measure. We assume that the uncertainty in our model arises from a vector of independent
% random variables $\xi=\left(\xi_{1},\dots,\xi_{m_{\xi}}\right)^{T}$ defined as $\xi:\Omega\rightarrow\Phi\subset\mathbb{R}^{m_{\xi}}$. Let~$\mathcal{B}(\Phi)$ be the Borel $\sigma$-algebra on $\Phi$ induced by $\xi$, and let $\mu$ denote the corresponding probability measure. The expectation of the product of measurable functions~$u$ and $v$ depending on $\xi$ defines a Hilbert space $L^{2}\left(\Phi\right)\coloneqq L^{2}\left(\Phi,\mathcal{B}(\Phi),\mu\right)$ equipped with the inner product
% $
% \left\langle u,v\right\rangle =\mathbb{E}\left[  uv\right]  =\int_{\Phi
% }u\left(  \xi\right)  v\left(  \xi\right)  \,d\mu\left(  \xi\right)  ,
% $
% where $\mathbb{E}$ denotes the mathematical expectation.
Following the standard framework~\cite{Ghanem-1991-SFE, Xiu-2010-NMS}, let $(\Omega, \mathcal{F}, \mathcal{P})$ be a complete probability space with $\xi: \Omega \to \Phi \subset \mathbb{R}^{m_\xi}$ a vector of independent random variables. We consider the Hilbert space $L^2(\Phi)$ equipped with the inner product $\langle u,v \rangle = \mathbb{E}[uv]$, where $\mathbb{E}$ denotes the mathematical expectation with respect to the probability measure $\mu$ induced by $\xi$.

We consider the
steady-state optimal control problem given by
\begin{equation}
    \min_{y,u}{J}(y,u):=\,\frac{1}{2}\int_{\Phi}\int_{\mathcal{D}} |y - y_{d}|^2 dx\,d\mu(\xi)
    +\frac{\beta}{2}\int_{\Phi}\int_{\mathcal{D}} |u|^2 \,dx\,d\mu(\xi)+
    % Var
    \frac{\gamma}{2} \int_{\mathcal{D}}|\sigma(y)|^2 \,dx,
    \label{eq:problem-stoch}
\end{equation}
subject~to
\begin{equation}
    \begin{cases}
        -\nabla\cdot\left(\Bbbk\left(x,\xi\right)\nabla y\left(x,\xi\right)\right)  = u\left(x,\xi\right), & \quad\text{in }\mathcal{D}\times \Phi,\label{eq:pde-stoch} \\
        y\left(x,\xi\right) = g\left(x\right),                                                             & \quad\text{on }\partial \mathcal{D}\times \Phi,
    \end{cases}
\end{equation}
where $y$ is the state, $y_{d}$ is the target state, and $u$ is the (distributed) control.
The parameter~$\gamma$ penalizes the variance~$\sigma^2(y)$ of the state $y$.
% where $$% Standard deviation (scalar random variable g)
% \sigma(y)= \left( \int_{\Omega} \bigl(y - \mathbb{E}[y]\bigr)^2 \, d\mathbb{P}(\omega) \right)^{1/2}.
% $$
Observe that both the state~$y$ and control~$u$ are stochastic. In view of the Doob-Dynkin lemma (see, e.g.,~\cite{Babuska-2004-GFE}),
both~$y$ and~$u$ admit the same parametric dependence on~$\xi$.

%
% Reordered Section 2 with Affine Expansion
%
In computations, we work with a finite dimensional subspace $\mathcal{T}_{p}\mathcal{\subset}L^{2}\left(  \Phi,\mathcal{B}(\Phi),\mu\right)$,
spanned by a set of generalized polynomial chaos (gPC) functions~$\left\{  \psi_{\ell}\left(  \xi\right)  \right\}_{\ell=1}^{n_A}$ with
\begin{equation}
    \langle \psi_{\ell},\psi_{k}\rangle=\mathbb{E}\left[  \psi_{\ell}\psi_{k}\right]=\int_\Phi \psi_\ell(\xi)\psi_k(\xi)d\mu=\delta_{\ell k},
\end{equation}
where $\left\{  \psi_{i}\left(  \xi\right)  \right\}  _{i=1}^{n_A}$ is a set of $m_\xi-$dimensional, $p$-order Hermite polynomials.~$\psi_1(\xi)=1$, and $\mathbb{E}[\psi_k(\xi)] = 0$ for $k > 1$.

It is assumed that~$\Bbbk(x,\xi)$ is bounded away from zero and infinity, i.e.,~$0 < \Bbbk_{\min} \le \Bbbk(x,\xi) \le \Bbbk_{\max} < \infty$ for some constants~$\Bbbk_{\min}$, $\Bbbk_{\max}$.
We assume the stochastic PDE coefficient~$\Bbbk(x,\xi)$ is represented by a gPC expansion
\begin{equation}\label{KLE}
    \Bbbk(x,\xi) = \sum_{i=1}^{n_A}\kappa_i(x)\psi_i(\xi), %\quad %\text{where } \kappa_i(x) = \frac{\mathbb{E}[\Bbbk(x)\psi_i(\xi)]}{\mathbb{E}[\psi_i^2]},
\end{equation}
Similarly, applying the stochastic Galerkin method to (\ref{eq:problem-stoch})--(\ref{eq:pde-stoch}),
both the state~$y$ and control~$u$ are expanded as
\begin{equation}
    v=\sum\nolimits_{i=1}^{n_h}\sum\nolimits_{k=1}^{n_\xi}v_{ik}\phi_{i}(x)\psi_{k}(\xi), \qquad v = y \, \text{or} \, u.
    \label{eq:gPC-yu}%
\end{equation}

Let $f(x,\xi)=f_0(x)+\sum_{i=1}^{n_{\xi}}\sqrt{\theta_i}f_i(x)\xi_i$ be a truncated Karhunen-Lo\`{e}ve(KL) expansion of a Gaussian process defined on $\mathcal{D}$, where $\xi_i$ are independent, identically distributed random variables, $f_0(x)$ is the mean function, and $(\theta_i,f_i(x))$ is the $i$-th eigenpair (eigenvalue and eigenfunctions, respectively) of the covariance function $C_f(x,y)$ such that
\begin{equation}
    \label{Czxy}
    \int_\mathcal{D} C_f(x,y)\varphi_i(y)dy=\theta_i\varphi_i(x),
\end{equation}
where $f_i(x)=\sqrt{\theta_i}\varphi_i(x)$. In particular, if  $\Bbbk(x,\xi)=\exp{[f(x,\xi)]}$, then following~\cite{Ghanem-1999-NGS}
\[
    \kappa_i\zhl{(x)}=\frac{\mathbb{E}[\psi_i(\xi-f)]}{\mathbb{E}[\psi_i^2]}\exp{\left[ f_0(x)+\frac{1}{2}\sum_{j=1}^{m_\xi}f^2_j(x) \right]}.
\]
According to~\cite{Matthies-2005-GML}, to guarantee a complete representation of the lognormal random field, the variables $y$ and $u$ are defined in the space $\mathcal{T}_{p}$ with dimension $n_\xi=\binom{m_{\xi}+p}{p}$,
% Consequently, the order of polynomial expansion of $\Bbbk(x,\xi)$ in (4) should be twice the order of the expansion of the solutionu(x,ξ)in (5)
%where $p$ is the total degree of polynomials and $m_\xi$ is the number of random variables. 
and the dimension of $\Bbbk(x,\xi)$ is given by $n_A=\binom{m_{\xi}+2p}{2p}$.
% \label{sec:stochastic-model} 

% After the application of the finite element method in the physical space to the variational form of (\ref{eq:pde-stoch}), the expansion~\eqref{KLE} leads to a stochastic stiffness matrix 
% \begin{equation}
%     \label{eq:affine-A}
%     A(\xi) = \sum_{i=1}^{n_A} A_i \psi_i(\xi),
% \end{equation}
% where $A_i$ are stiffness matrices corresponding to the coefficient modes $\kappa_i(x)$.

%where $v_{ik}$ represents the deterministic vector of coefficients or degrees of freedom~(dof),  
%\[
%v_{1:n_h,1},v_{1:n_h,2},\dots,v_{1:n_h,n_\xi}.
%\]
Using the stochastic Galerkin finite element framework to discretize problem~(\ref{eq:problem-stoch})--(\ref{eq:pde-stoch}), we get
\[
    \min_{\mathbf{y},\mathbf{u}}J\left(\mathbf{y},\mathbf{u}\right)  =\frac{1}{2}\left(\mathbf{y}-\mathbf{y}_{d}\right)^{T}\mathcal{M}\left( \mathbf{y}-\mathbf{y}_{d}\right)  +\frac{\gamma
    }{2}\mathbf{y}^{T}\mathcal{M}_{\sigma}\mathbf{y}+\frac{\beta}{2}\mathbf{u}%
    ^{T}\mathcal{M}\mathbf{u},
\]
subject to
\begin{equation}
    -\mathcal{A} \mathbf{y}+\mathcal{M} \mathbf{u}=\mathbf{g}.
\end{equation}
Via the Lagrangian formulation, we can derive the first-order optimality conditions %using $H^{\gamma}%
%=H_{1}+\gamma H^{\sigma}$,~and
% \[
% %\begin{aligned}
% %    \mathcal{M}_{\gamma}&=\mathcal{M+\gamma M}_{\sigma}=H^{\gamma}\otimes M,\\
% %    \mathcal{A}&=\sum_{\ell=1}^{n_{A}}A_{\ell}\otimes\psi_{\ell}(\xi)=\sum_{\ell=1}^{n_{A}}H_{\ell}\otimes A_{\ell},
% %\end{aligned}
%     \mathcal{M}_{\gamma}=H^{\gamma}\otimes M, \qquad \mathcal{A}=\sum\nolimits_{\ell=1}^{n_{A}}H_{\ell}\otimes A_{\ell},
% \]
% we get the KKT system of equations,  
%\begin{align}
%\mathcal{M}_{\gamma}\mathbf{y}-\mathcal{A}^{T}\lambda &  =\mathcal{M}\mathbf{y}_{d},\label{eq:stoch-adjoint}\\
%\beta\mathcal{M}\mathbf{u}+\mathcal{M}^{T}\lambda&=0,\label{eq:stoch-gradient}\\
%-\mathcal{A}\mathbf{y}+\mathcal{M}\mathbf{u} &  =\mathbf{g},\label{eq:stoch-state}%
%\end{align}
%and 
% which can be written in the matrix form as 
\begin{equation}
    \left[
        \begin{array}
            [c]{ccc}%
            \mathcal{M}_{\gamma} & 0                & -\mathcal{A}^{T} \\
            0                    & \beta\mathcal{M} & \mathcal{M}^{T}  \\
            -\mathcal{A}         & \mathcal{M}      & 0
        \end{array}
        \right]  \left[
        \begin{array}
            [c]{c}%
            \mathbf{y} \\
            \mathbf{u} \\
            \zhl{\boldsymbol{\lambda}}%
        \end{array}
        \right]  =\left[
        \begin{array}
            [c]{c}%
            \mathcal{M}\mathbf{y}_{d} \\
            0                         \\
            \mathbf{g}%
        \end{array}
        \right],
    \label{eq:stoch-system}%
\end{equation}
where $\boldsymbol{\lambda}$ is a Lagrangian multiplier, $\mathcal{M}_{\gamma}=\mathcal{M+\gamma M}_{\sigma}=H^{\gamma}\otimes M$,  $H^{\gamma}=H_{1}+\gamma H^{\sigma}$, and similarly, $h^\gamma$ is defined as the entry of $H^\gamma$.  In particular, 
\begin{equation}
    %\begin{aligned}\label{eq:def-Mcal}
    %        \mathcal{M}  &  =H_{1}\otimes M,\\
    %        \mathcal{M}_{\sigma}  &  =H^{\sigma}\otimes M.
    %\end{aligned}   
    \mathcal{M}=H_{1}\otimes M, \qquad \mathcal{M}_{\sigma}=H^{\sigma}\otimes M,
    \label{eq:def-Mcal}
\end{equation}
where 
%Using the expansions (\ref{KLE}) and (\ref{eq:gPC-yu})  yields the global stiffness matrix system
the mass matrix, $M$ is given by $$M_{a,b}=\int_D\phi_a(x)\phi_b(x) dx,  \quad a,b=1,2,\dots,n_h,$$ 
for a suitably chosen finite element basis, $\{\phi_i(x)\}_{i=1}^{n_h}.$
%The stochastic Galerkin matrix in~(\ref{eq:stoch-system}) is symmetric, indefinite, and in general very large. 
%A preconditioner for~(\ref{eq:stoch-system}), proposed by Benner et al.~\cite{Benner-2016-BDP}, was formulated as
 The global stiffness matrix $\mathcal{A}$ is given by $\mathcal{A}=\sum_{\ell=1}^{n_A}H_\ell\otimes A_\ell$.
%, where the matrices $H_\ell$ are defined below.
Here, the matrices $H_\ell$ are defined by their entries$(H_\ell)_{jk} = h_{\ell j k}$, with $h_{ijk} = \mathbb{E}[\psi_i \psi_j \psi_k]$,
where we note that all matrices$~H_{\ell}$ are symmetric, and%
\[
H^{\sigma}=\operatorname{diag}\left(  0,h_{1,jj}\right)  ,\quad j=2,\dots,n_{\xi},
\]
which is obtained from~$H_{1}$\ by setting~$h_{1,11}=0$. 
% This system can be written as $\mathcal{A} \bar{y} = \bar{F}$, where $\mathcal{A}=\sum_{\ell=1}^{n_A} H_\ell\otimes A_\ell$, and the matrices $H_\ell$ are defined by $(H_\ell)_{jk} = c_{\ell j k}$.
% We note that all matrices $A_\ell$ share the same sparsity pattern. The resulting linear system can be expressed as $\mathcal{A} \bar{y} = \bar{F}$, where the global stiffness matrix is given by $\mathcal{A}=\sum_{\ell=1}^{n_A} H_\ell\otimes A_\ell$. Here, the matrices $H_\ell$ are defined by their entries $(H_\ell)_{jk} = c_{\ell j k}$, with $c_{ijk} = \mathbb{E}[\psi_i \psi_j \psi_k]$, and $\bar{F}$ represents the projection of the right-hand side.
% \begin{equation}
% \label{eq:gPC-A}
% % {\color{black}A_\ell=[(A_\ell)_{ab}],\quad (A_\ell)_{ab}=\int_{\mathcal{D}}\zhl{\kappa_\ell(x)}\nabla \phi_a(x)\cdot\nabla \phi_b(x) \,dx, \quad a,b=1,2,\dots,n_h.}
% \sum_{j=1}^{n_\xi} \sum_{i=1}^{n_A} c_{ijk} A_i y_j = F_k, \quad k=1,\dots,n_\xi,
% \end{equation}
%
Note  also that in our settings~$H_{1}=I_{n_\xi}$, i.e., it is an identity matrix of size~$n_\xi$.
Moreover, all matrices $A_\ell$ share the same sparsity pattern and are given by
% Substituting the expansions (\ref{KLE}) and (\ref{eq:gPC-yu}) into the semi-discretized PDE yields the global stiffness matrix system
\begin{equation}
    {\color{black}A_\ell=[(A_\ell)_{ab}],\quad (A_\ell)_{ab}=\int_{\mathcal{D}}\zhl{\kappa_\ell(x)}\nabla \phi_a(x)\cdot\nabla \phi_b(x) \,dx, \quad a,b=1,2,\dots,n_h.}
    \label{eq:gPC-A}%
\end{equation}

% In implementation, we use the \emph{matricized} format, which utilizes the isomorphism between $\mathbb{R}^{n_h n_\xi}$ and $\mathbb{R}^{n_h \times n_\xi}$, 
% defined via the operators $\operatorname{vec}$ and $\operatorname{mat}$. 
% Specifically, 
% \begin{equation}
% \bar{V}=\operatorname{mat}(\bar{v})=\left[  v_{1},v_{2},\ldots,v_{n_\xi%
% }\right]  \in%
% %TCIMACRO{\U{211d} }%
% %BeginExpansion
% \mathbb{R}
% %EndExpansion
% ^{n_h\times n_\xi}, \label{eq:m}%
% \end{equation}
% where the column~$k$ contains the coefficients associated with the basis
% function$~\psi_{k}$, and~$\bar{v} = \operatorname{vec}\left(\bar{V}\right) \in \mathbb{R}^{n_h n_\xi}$.

% We will also use the notation
% \[
% H_{\ell}=\left[  h_{\ell,jk}\right]  ,\quad h_{\ell,jk}\equiv\mathbb{E}\left[
% \psi_{\ell}\psi_{j}\psi_{k}\right]  ,\quad\ell=1,\dots,n_{A},\;j,k=1,\dots
% ,n_\xi,
% \]
% where we note that all matrices$~H_{\ell}$ are symmetric, and%
%which is obtained from~$H_{1}$\ by setting~$h_{1,11}=0$. 

We note that the coefficient matrix in~(\ref{eq:stoch-system}) is symmetric, indefinite, and in general very large. It is ill-conditioned, and therefore, a good preconditioner is required to solve the system efficiently. Next, we introduce a preconditioner to tackle this problem.

\subsection{Schur complement-based preconditioner}
A block-diagonal preconditioner for~(\ref{eq:stoch-system}) (see, e.g. Benner et al.~\cite{Benner-2016-BDP}) is given by
%\mathcal{P}_\text{hGSoc}
\begin{equation}
    %\mathcal{P}=
    \mathcal{P}:=
    \begin{bmatrix}
        \mathcal{M}_{\gamma} & 0                & 0                        \\
        0                    & \beta\mathcal{M} & 0                        \\
        0                    & 0                & \mathcal{S}_\text{exact} %
    \end{bmatrix},
    \label{eq:stoch-prec}
\end{equation}
%where $\mathcal{S}_{\text{exact}}=\mathcal{A} \mathcal{M}_{\gamma}^{-1} \mathcal{A}^T+\frac{1}{\beta}\mathcal{M}$, but here we use another approximation $\mathcal{S}$ since it is quite difficult to invert.
where $\mathcal{S}_{\text{exact}}$ is the exact Schur complement
%. The exact Schur complement for this system is given by
\begin{equation}
    \label{S_exact}
    \mathcal{S}_\text{exact} = \mathcal{A} \mathcal{M}_{\gamma}^{-1} \mathcal{A}^T + \frac{1}{\beta}\mathcal{M}.
\end{equation}
The first two blocks corresponding to (scaling of) the mass matrix are block diagonal,
and the inverses are approximated by Chebyshev semi-iteration. %listed as Algorithm~\ref{alg:chebmass}. 
However, forming and applying the inverse of $\mathcal{S}_\text{exact}$ is computationally prohibitive. The primary difficulty stems from its additive structure. Therefore, the key to an efficient solution lies in designing an approximation $\mathcal{S}$ that is %both high quality and structurally convenient. 
spectrally equivalent to $\mathcal{S}_\text{exact}$ and easy to invert.
Following the approach in~\cite{Benner-2016-BDP}, we employ the approximation
\begin{equation}
    \label{eq:mathcalZ}
    \mathcal{S}=\mathcal{ZM}_{\gamma}^{-1}\mathcal{Z}^{T},\qquad\mathcal{Z}%
    =\mathcal{A}+\sqrt{\frac{1+\gamma}{\beta}}\mathcal{M}=\sum_{\ell=1}^{n_{A}%
        }H_{\ell}\otimes\tilde{A}_{\ell},
\end{equation}
where
\[
    \tilde{A}_{1}=A_{1}+\sqrt{\frac{1+\gamma}{\beta}}M,\qquad\tilde{A}_{\ell
    }=A_{\ell},\quad\ell=2,\dots,n_{A}.
\]
Since $\mathcal{Z}$ is symmetric, we have~$\mathcal{S}^{-1}=\mathcal{Z}^{-1}\mathcal{M}_\gamma\mathcal{Z}^{-1}$.
% An approximation of $\mathcal{Z}^{-1}$ by the mean-based preconditioner~\cite{Pellissetti-2000-ISS,Powell-2009-BDP} and its practical variants  were studied 
In \cite{Benner-2016-BDP}, the authors studied a mean-based preconditioner derived from  (\ref{eq:mathcalZ}) by dropping all the terms of $\mathcal{Z}$ except the first term; that is, $\mathcal{Z}\approx H_1 \otimes \tilde{A}_{1}$. However, this mean-based preconditioner is less effective when the standard deviation of the uncertain parameters increases.

In this study, we overcome this shortcoming by hierarchical preconditioning introduced in~\cite{Sousedik-2014-THP,Sousedik-2014-HSC}  in the context of forward problems. More specifically, unlike ~\cite{Sousedik-2014-THP,Sousedik-2014-HSC}, we extend this strategy to the stochastic optimal control problem,  and we also study the spectral properties of the new preconditioner.

To that end,  observe first that, by
applying the identity
\begin{equation}
    \left(  V\otimes W\right)  \operatorname{vec}(X)=\operatorname{vec}%
    (WXV^{T}),\label{eq:mat-vec}%
\end{equation}
to~(\ref{eq:stoch-system}) yields the matricized system~\cite[eq.~(53)]{Benner-2016-BDP}:
\begin{equation}
    \begin{aligned} M\bar{Y}H^{\gamma}-\sum_{\ell=1}^{n_{A}}\tilde{A}_{\ell}^{T}\bar{\Lambda}H_{\ell} &=M\bar {Y}_{d}H_{1} , \\ \beta M\bar{U}H_{1}+M^{T}\bar{\Lambda}H_{1} &=0 , \\ -\sum_{\ell=1}^{n_{A}}\tilde{A}_{\ell}^{T}\bar{Y}H_{\ell}+M\bar{U}H_{1} &=\bar{G} . \end{aligned}\label{eq:stoch-system-m}%
\end{equation}
In addition, the matrix-vector of $\mathcal{Z}\bar{v}$ in~(\ref{eq:mathcalZ}) reads
\begin{equation}
    \mathcal{Z}\bar{v}=\operatorname{vec}\left(\sum\nolimits_{\ell=1}^{n_{A}}\tilde{A}_{\ell
    }\bar{V}H_{\ell}\right).\label{eq:Zv}%
\end{equation}

\subsection{Hierarchical Gauss-Seidel preconditioner}
\label{sec:preconditioning}

We first recall the preconditioner for the forward PDE problem from~\cite{Sousedik-2014-THP}. However, in this work,  we present it in the matricized format as Algorithm~\ref{alg:hGS}.
% Moreover, we formulate the corresponding preconditioner for the optimal control problem as Algorithm~\ref{alg:hGSoc-1}--\ref{alg:hGSoc-2}. In implementation, we use the \emph{matricized} format, utilizing the isomorphism between $\mathbb{R}^{n_h n_\xi}$ and $\mathbb{R}^{n_h \times n_\xi}$. Specifically, $\bar{V}=\operatorname{mat}(\bar{v})=\left[ v_{1},\dots,v_{n_\xi}\right]\in \mathbb{R}^{n_h\times n_\xi}$, where~$v_k$ contains coefficients associated with $\psi_k$, and $\bar{v} = \operatorname{vec}(\bar{V})$.
To set the notation,
 we note that the \emph{matricized} format, which utilizes the isomorphism between $\mathbb{R}^{n_h n_\xi}$ and $\mathbb{R}^{n_h \times n_\xi}$, 
defined via the operators $\operatorname{vec}$ and $\operatorname{mat}$. 
Specifically, 
\begin{equation}
\bar{V}=\operatorname{mat}(\bar{v})=\left[  v_{1},v_{2},\ldots,v_{n_{\xi}%
}\right]  \in%
%TCIMACRO{\U{211d} }%
%BeginExpansion
\mathbb{R}
%EndExpansion
^{n_h\times n_{\xi}}, \label{eq:m}%
\end{equation}
where the column~$k$ contains the coefficients associated with the basis
function$~\psi_{k}$, and
$\bar{v} = \operatorname{vec}\left(\bar{V}\right) \in \mathbb{R}^{n_h n_\xi}$.
We will use lowercase letters for the \emph{vectorized} representation and uppercase letters for the \emph{matricized} counterpart; 
so, e.g., $\bar{R} = \operatorname{mat}(\bar{r})$, etc.

Moreover, we will denote by$~\bar{V}_{\left(  i:n\right) }$ a submatrix of~$\bar{V}$
containing columns~$i,i+1,\dots,n$, and, in particular,
$\bar{V}=\bar{V}_{\left(  1:n_\xi\right) }$. There are two components of the
preconditioner. The first component consists of block-diagonal solves with
blocks of varying sizes. The second component is used in the setup of the
right-hand sides for the solves, and consists of matrix-vector products by
certain subblocks of the stochastic Galerkin matrix by vectors of
corresponding sizes. To this end, we will write~$\left[  h_{\tau,(\ell)(k)}\right]$,
with~$(\ell)$\ and $(k)$ denoting a set of (consecutive) rows
and columns of matrix~$H_{\tau}$ so that, in particular,~${H}_{\tau}=\left[h_{\tau,(1:n_\xi)(1:n_\xi)}\right]$.
Let us also denote $\bar{v}_{(\ell)}=\operatorname{vec}(\bar{V}_{(\ell)})$.
Then, the matrix-vector products can be written, cf.~(\ref{eq:mat-vec}) and noting the symmetry of$~H_{\tau}$, as
\begin{equation}
    \bar{v}_{(\ell)}=\sum\nolimits_{\tau\in\mathcal{I}_{\tau}}(\left[  h_{\tau,(\ell)(k)}\right]
    \otimes\tilde{A}_{\tau})\bar{u}_{(k)}\quad\Leftrightarrow\quad\bar{V}_{(\ell)}
    =\sum\nolimits_{\tau\in\mathcal{I}_{\tau}}\tilde{A}_{\tau}\bar{U}_{(k)}\left[  h_{\tau,(k)(\ell)}\right]  , \label{eq:matvec-block}%
\end{equation}
where$~\mathcal{I}_{\tau}$\ is an index set$~\mathcal{I}_{\tau}\subseteq\left\{
    1,\dots,n_{\tau}\right\}  $ indicating that the matrix-vector products may be
truncated. Possible strategies for truncation are discussed
in~\cite{Sousedik-2014-THP}. In this study, we use $\mathcal{I}_{\tau}=\left\{
    1,\dots,n_{\tau}\right\}  $ with $n_{\tau}=\binom{m_{\xi}+p_{\tau}}{p_{\tau}}$ for some
$p_{\tau}\leq p$. In particular, we set $\tau=\{0,1,2\}$ and with $\mathcal{I}_{\tau}=\emptyset$ both preconditioners in Algorithms~\ref{alg:hGS}
and~\ref{alg:hGSoc-1}--\ref{alg:hGSoc-2} reduce to mean-based variants.
We also note that, since the initial guess is zero in Algorithm~\ref{alg:hGS}, the
multiplications by$~\mathcal{F}_{1}$\ and $\mathcal{F}_{d+1}$ vanish
from~(\ref{eq:alg-hGS1})--(\ref{eq:alg-hGS2}).

\begin{algorithm}[t]
    \caption{{\cite[Algorithm~3]{Sousedik-2014-THP}} Hierarchical Gauss-Seidel preconditioner (hGS)}
    \label{alg:hGS}
    The preconditioner $\mathcal{Z}_{hGS}:\bar{R}\longmapsto\bar{V}$ is defined as follows.
    \begin{algorithmic}[1]
        \State Set the initial solution$~\bar{V}$ to zero and update in the following steps:
        \State Solve
        \begin{equation}
            \tilde{A}_1\bar{V}_{\left(  1\right)  }=\bar{R}_{\left(  1\right)  }-\mathcal{F}%
            _{1},\qquad\text{where }\mathcal{F}_{1}=\sum_{\tau\in\mathcal{I}_{\tau}}\tilde{A}_{\tau}\bar
            {V}_{\left(  2:n_\xi\right)  }\left[  h_{\tau,\left(  2:n_\xi\right)  \left(
                        1\right)  }\right]  .\label{eq:alg-hGS1}
        \end{equation}
        \For{$d=1,\ldots, p-1$}
        \State Set
        $\ell =\left( n_{\ell }+1:n_{u}\right) ,\text{ where }n_{\ell }=\binom{m_{\xi
                    }+d-1}{d-1}\text{ and }n_{u}=\binom{m_{\xi }+d}{d}$.
        \State Solve
        \begin{equation}
            \tilde{A}_1\bar{V}_{\left(  \ell\right)  }=\bar{R}_{\left(  \ell\right)
            }-\mathcal{E}_{d+1}-\mathcal{F}_{d+1},\label{eq:alg-hGS2}%
        \end{equation}
        where%
        \begin{equation*}
            \mathcal{E}_{d+1}   =\sum_{\tau\in\mathcal{I}_{\tau}}\tilde{A}_{\tau}\bar{V}_{\left(
                    1:n_{\ell}\right)  }\left[  h_{\tau,\left(  1:n_{\ell}\right)  \left(
                        \ell\right)  }\right]  , \quad
            \mathcal{F}_{d+1}   =\sum_{\tau\in\mathcal{I}_{\tau}}\tilde{A}_{\tau}\bar{V}_{\left(
                    n_{u}+1:n_\xi\right)  }\left[  h_{\tau,\left(  n_{u}+1:n_\xi\right)  \left(
                        \ell\right)  }\right]  .
        \end{equation*}
        \EndFor
        \State Set $\ell =\left( n_{u}+1:n_{\xi }\right)$.
        \State Solve
        \begin{equation*}
            \tilde{A}_1\bar{V}_{\left(  \ell\right)  }=\bar{R}_{\left(  \ell\right)
            }-\mathcal{E}_{p+1},\qquad\text{where }\mathcal{E}_{p+1}=\sum_{t\in
                \mathcal{I}_{\tau}}\tilde{A}_{\tau}\bar{V}_{\left(  1:n_{u}\right)  }\left[  h_{\tau,\left(
                        1:n_{u}\right)  \left(  \ell\right)  }\right]  ,
        \end{equation*}
        \For{$d=p-1,\ldots, 1$}
        \State Set $\ell =\left( n_{\ell }+1:n_{u}\right) ,\text{ where }n_{\ell }=\binom{m_{\xi
                    }+d-1}{d-1}\text{ and }n_{u}=\binom{m_{\xi }+d}{d}$.
        \State Solve~(\ref{eq:alg-hGS2}).
        \EndFor
        \State Solve~(\ref{eq:alg-hGS1}).
    \end{algorithmic}
\end{algorithm}

Next, we apply the hGS strategy to form a preconditioner for the KKT system.
The preconditioner is formulated as Algorithm~\ref{alg:hGSoc-1}--\ref{alg:hGSoc-2}.
It adapts the core idea of Algorithm~\ref{alg:hGS} to handle the coupled variables corresponding
to the state~$\bar{V}^{Y}$, control~$\bar{V}^{U}$, and adjoint~$\bar{V}^{\Lambda}$ simultaneously at each hierarchical level.
\zh{
A key computational step of Algorithm~\ref{alg:hGSoc-1}--\ref{alg:hGSoc-2} entails a solve with $\tilde{P}=\operatorname{blkdiag}\begin{bmatrix}
M,\beta M,\tilde{A}_{1}\end{bmatrix}$, which serves as an auxiliary block-diagonal preconditioner derived from the blocks of the deterministic KKT system.}
Since this system is relatively small and constant across all hierarchical levels,
it can be handled efficiently, for instance, by computing a direct factorization of $\tilde{P}$ once and reusing it for all subsequent solves.
The overall performance of the %hGSoc 
preconditioner is thus determined by the cost of these deterministic solves
and the number of truncated off-diagonal matrix-vector products.

\begin{algorithm}[hptb]
    \caption{hGS preconditioner for the optimal control problem (hGSoc)}
    \label{alg:hGSoc-1}
    The preconditioner $\mathcal{P}_{hGSoc}:\left(  \bar{R}^{Y},\bar{R}^{U}%
        ,\bar{R}^{\Lambda}\right)  \longmapsto\left(  \bar{V}^{Y},\bar{V}^{U},\bar
        {V}^{\Lambda}\right)  $ is defined as follows.
    \begin{algorithmic}[1]
        \State Set the initial solution $\left(  \bar{V}^{Y},\bar{V}^{U},\bar{V}^{\Lambda
            }\right) $ to zero and update in the following steps:
        \State Solve
        \begin{equation}
            \tilde{P}\left[
                \begin{array}
                    [c]{c}%
                    \bar{V}_{\left(  1\right)  }^{Y} \\
                    \bar{V}_{\left(  1\right)  }^{U} \\
                    \bar{V}_{\left(  1\right)  }^{\Lambda}%
                \end{array}
                \right]  =\left[
                \begin{array}
                    [c]{c}%
                    \bar{R}_{\left(  1\right)  }^{Y}-\mathcal{C}_{1}+\mathcal{D}_{1} \\
                    \bar{R}_{\left(  1\right)  }^{U}-\mathcal{E}_{1}-\mathcal{F}_{1} \\
                    \bar{R}_{\left(  1\right)  }^{\Lambda}+\mathcal{G}_{1}-\mathcal{H}_{1}%
                \end{array}
                \right]  ,\label{eq:alg-hGSoc1}%
        \end{equation}
        where%
        \begin{align*}
            \mathcal{C}_{1}                                                  & =M\bar{V}_{\left(2:n_\xi\right)  }^{Y}\left[
            h_{\left(  2:n_\xi\right)  \left(  1\right)  }^{\gamma}\right] , &
            \mathcal{D}_{1}                                                  & =\sum_{\tau\in\mathcal{I}_{\tau}}\tilde{A}_{\tau}\bar{V}_{\left(  2:n_{\xi
                    }\right)  }^{\Lambda}\left[  h_{\tau,\left(  2:n_\xi\right)  \left(  1\right)
            }\right] ,                                                                                                                                                                   \\
            \mathcal{E}_{1}                                                  & =\beta M\bar{V}_{\left(  2:n_\xi\right)  }^{U}\left[
            h_{1,\left(  2:n_\xi\right)  \left(  1\right)  }\right] ,        & \mathcal{F}_{1}
            & =M\bar{V}_{\left(  2:n_\xi\right)  }^{\Lambda}\left[  h_{1,\left(2:n_\xi\right)  \left(  1\right)  }\right] ,                                                                                                                                 \\
            \mathcal{G}_{1}                                                  & =\sum_{\tau\in\mathcal{I}_{\tau}}\tilde{A}_{\tau}\bar{V}_{\left(  2:n_{\xi
                    }\right)  }^{Y}\left[  h_{\tau,\left(  2:n_\xi\right)  \left(  1\right)
            }\right] ,                                                       & \mathcal{H}_{1}                                                    & =M\bar{V}_{\left(  2:n_\xi\right)  } %
            ^{U}\left[  h_{1,\left(  2:n_\xi\right)  \left(  1\right)  }\right] .
        \end{align*}
        \For{$d=1,\ldots, p-1$}
        \State Set
        $\ell =\left( n_{\ell }+1:n_{u}\right) ,\text{ where }n_{\ell }=\binom{m_{\xi
                    }+d-1}{d-1}\text{ and }n_{u}=\binom{m_{\xi }+d}{d}$.
        \State Solve
        \begin{equation}
            \tilde{P}\left[
                \begin{array}
                    [c]{c}%
                    \bar{V}_{\left(  \ell\right)  }^{Y} \\
                    \bar{V}_{\left(  \ell\right)  }^{U} \\
                    \bar{V}_{\left(  \ell\right)  }^{\Lambda}%
                \end{array}
                \right]  =\left[
                \begin{array}
                    [c]{c}%
                    \bar{R}_{\left(  \ell\right)  }^{Y}-\mathcal{C}_{d+1}+\mathcal{D}_{d+1} \\
                    \bar{R}_{\left(  \ell\right)  }^{U}-\mathcal{E}_{d+1}-\mathcal{F}_{d+1} \\
                    \bar{R}_{\left(  \ell\right)  }^{\Lambda}+\mathcal{G}_{d+1}-\mathcal{H}_{d+1}%
                \end{array}
                \right]  ,\label{eq:alg-hGSoc2}%
        \end{equation}
        where
        \begin{align*}
            \mathcal{C}_{d+1} & =M\left(  \bar{V}_{\left(  1:n_{\ell}\right)  }                  %
            ^{Y}\left[  h_{\left(  1:n_{\ell}\right)  \left(  \ell\right)  }^{\gamma}\right]  +\bar{V}_{\left(  n_{u}+1:n_\xi\right)  }^{Y}\left[  h_{\left(
            n_{u}+1:n_\xi\right)  \left(  \ell\right)  }^{\gamma}\right]  \right) ,              \\
            \mathcal{D}_{d+1} & =\sum_{\tau\in\mathcal{I}_{\tau}}\tilde{A}_{\tau}\left(  \bar{V}_{\left(
                1:n_{\ell}\right)  }^{\Lambda}\left[  h_{\tau,\left(  1:n_{\ell}\right)  \left(
                    \ell\right)  }\right]  +\bar{V}_{\left(  n_{u}+1:n_\xi\right)  }^{\Lambda
            }\left[  h_{\tau,\left(  n_{u}+1:n_\xi\right)  \left(  \ell\right)  }\right]
            \right) ,                                                                            \\
            \mathcal{E}_{d+1} & =\beta M\left(  \bar{V}_{\left(  1:n_{\ell}\right)
            }^{U}\left[  h_{1,\left(  1:n_{\ell}\right)  \left(  \ell\right)  }\right]
            +\bar{V}_{\left(  n_{u}+1:n_\xi\right)  }^{U}\left[  h_{1,\left(
            n_{u}+1:n_\xi\right)  \left(  \ell\right)  }\right]  \right) ,                       \\
            \mathcal{F}_{d+1} & =M\left(  \bar{V}_{\left(  1:n_{\ell}\right)
            }^{\Lambda}\left[  h_{1,\left(  1:n_{\ell}\right)  \left(  \ell\right)
                }\right]  +\bar{V}_{\left(  n_{u}+1:n_\xi\right)  }^{\Lambda}\left[
                h_{1,\left(  n_{u}+1:n_\xi\right)  \left(  \ell\right)  }\right]  \right) ,
            \\
            \mathcal{G}_{d+1} & =\sum_{\tau\in\mathcal{I}_{\tau}}\tilde{A}_{\tau}\left(  \bar{V}_{\left(
                1:n_{\ell}\right)  }^{Y}\left[  h_{\tau,\left(  1:n_{\ell}\right)  \left(
                    \ell\right)  }\right]  +\bar{V}_{\left(  n_{u}+1:n_\xi\right)  }^{Y}\left[
                h_{\tau,\left(  n_{u}+1:n_\xi\right)  \left(  \ell\right)  }\right]  \right) ,
            \\
            \mathcal{H}_{d+1} & =M\left(  \bar{V}_{\left(  1:n_{\ell}\right)  }                  %
            ^{U}\left[  h_{1,\left(  1:n_{\ell}\right)  \left(  \ell\right)  }\right]
            +\bar{V}_{\left(  n_{u}+1:n_\xi\right)  }^{U}\left[  h_{1,\left(
                    n_{u}+1:n_\xi\right)  \left(  \ell\right)  }\right]  \right) .
        \end{align*}
        \EndFor
        \algstore{split_here}
    \end{algorithmic}
\end{algorithm}

\begin{algorithm}[hptb]
    \caption{hGS preconditioner for the optimal control problem (hGSoc), cont'd}
    \label{alg:hGSoc-2}
    \begin{algorithmic}[1]
        \algrestore{split_here}
        \State Set $\ell =\left( n_{u}+1:n_{\xi }\right)$.
        \State Solve
        \[
            \tilde{P}\left[
                \begin{array}
                    [c]{c}%
                    \bar{V}_{\left(  \ell\right)  }^{Y} \\
                    \bar{V}_{\left(  \ell\right)  }^{U} \\
                    \bar{V}_{\left(  \ell\right)  }^{\Lambda}%
                \end{array}
                \right]  =\left[
                \begin{array}
                    [c]{c}%
                    \bar{R}_{\left(  \ell\right)  }^{Y}-\mathcal{C}_{p+1}+\mathcal{D}_{p+1} \\
                    \bar{R}_{\left(  \ell\right)  }^{U}-\mathcal{E}_{p+1}-\mathcal{F}_{p+1} \\
                    \bar{R}_{\left(  \ell\right)  }^{\Lambda}+\mathcal{G}_{p+1}-\mathcal{H}_{p+1}%
                \end{array}
                \right]  ,
        \]
        where
        \begin{align*}
            \mathcal{C}_{p+1}                                                        & =M\bar{V}_{\left(  1:n_{u}\right)  }^{Y}\left[
            h_{\left(  1:n_{u}\right)  \left(  \ell\right)  }^{\gamma}\right] ,      &
            \mathcal{D}_{p+1}                                                        & =\sum_{\tau\in\mathcal{I}_{\tau}}A_{\tau}\bar{V}_{\left(
                    1:n_{u}\right)  }^{\Lambda}\left[  h_{\tau,\left(  1:n_{u}\right)  \left(
            \ell\right)  }\right] ,                                                                                                                                                    \\
            \mathcal{E}_{p+1}                                                        & =\beta M\bar{V}_{\left(  1:n_{u}\right)  }                                                      %
            ^{U}\left[  h_{1,\left(  1:n_{u}\right)  \left(  \ell\right)  }\right] , &
            \mathcal{F}_{p+1}                                                        & =M\bar{V}_{\left(  1:n_{u}\right)  }^{\Lambda
            }\left[  h_{1,\left(  1:n_{u}\right)  \left(  \ell\right)  }\right] ,                                                                                                      \\
            \mathcal{G}_{p+1}                                                        & =\sum_{\tau\in\mathcal{I}_{\tau}}A_{\tau}\bar{V}_{\left(
                    1:n_{u}\right)  }^{Y}\left[  h_{\tau,\left(  1:n_{u}\right)  \left(  \ell\right)
            }\right] ,                                                               & \mathcal{H}_{p+1}                                        & =M\bar{V}_{\left(  1:n_{u}\right)  } %
                ^{U}\left[  h_{1,\left(  1:n_{u}\right)  \left(  \ell\right)  }\right] .
        \end{align*}
        \For{$d=p-1,\ldots, 1$}
        \State Set $\ell =\left( n_{\ell }+1:n_{u}\right) ,\text{ where }n_{\ell }=\binom{m_{\xi
                    }+d-1}{d-1}\text{ and }n_{u}=\binom{m_{\xi }+d}{d}$.
        \State Solve~(\ref{eq:alg-hGSoc2}).
        \EndFor
        \State Solve~(\ref{eq:alg-hGSoc1}).
    \end{algorithmic}
\end{algorithm}

%\newpage
\section{Time-dependent problem}
\label{sec:Time-Dependent Problem}
\vspace{-.3cm}
%\subsection{Problem formulation}
\label{sec:stochastic-model-td}

The time-dependent optimal control problem is given by
%. We define the objective functional as
\begin{eqnarray}
    \label{eq:obj:stoch-time-dependent}
    %J\left(  y,u,t\right)  =
    \min_{y,u}\mathcal{J}(y,u)&=& \frac{1}{2}\int_0^T\int_\Phi\int_\mathcal{D} |y-y_{d}|^2 dx\,d\mu(\xi)\,dt+\frac{\beta}{2}\int_0^T\int_{\Phi}\int_\mathcal{D} |u|^2 dx\,d\mu(\xi)\,dt\nonumber\\
    && + \frac{\gamma}{2}\int_0^T\int_\Phi\int_\mathcal{D} |\sigma(y)|^2 dx\,d\mu(\xi)\,dt,
\end{eqnarray}
% \[
% J\left(  y,u,t\right)  =\frac{1}{2}\mathbb{E}\left[  \left\Vert y-y_{d}%
% \right\Vert _{L^{2}(D;0,T)\otimes L^2(\Omega)}\right]+\frac{\beta}%
% {2}\mathbb{E}\left[  \left\Vert u\right\Vert _{L^{2}(D;0,T)\otimes L^2(\Omega)}\right]  .
% \]
subject to %, $\mathbf{P}$-almost surely, to
\begin{equation}
    \label{eq:constraint:stochastic-time-dependent}
    \begin{cases}\begin{aligned}
            \frac{\partial y(t, \mathbf{x}, \xi)}{\partial t}-\nabla \cdot( \Bbbk(\mathbf{x}, \xi) \nabla y(t, \mathbf{x}, \xi)) & =u(t, \mathbf{x}, \xi) \text { in }(0, T] \times \mathcal{D} \times \Phi, \\
            y(t, \mathbf{x}, \xi)                                                                                                & =g \text { on }(0, T] \times \partial \mathcal{D} \times \Phi,            \\
            y(0, \mathbf{x}, \xi)                                                                                                & =y_0 \text { in } \mathcal{D} \times \Phi.
        \end{aligned}\end{cases}
\end{equation}
After the application of the stochastic Galerkin finite element discretization to (\ref{eq:obj:stoch-time-dependent}),
%--(\ref{eq:constraint:stochastic-time-dependent})
and using the trapezoidal rule for the time discretization, where $n_t=T/\tau$ is the number of time steps over the interval $[0,T]$
with time-step size $\tau$, we obtain
\begin{equation}\label{eq:obj:stoch-time-dependent-discretize}
    \min_{\mathbf{y},\mathbf{u}}   \mathcal{J}\left(  \mathbf{y},\mathbf{u}\right)=\frac{\tau}{2}(\mathbf{y}-\mathbf{y}_d)^T \left(D\otimes \mathcal{M}_{\gamma}\right)(\mathbf{y}-\mathbf{y}_d)+\frac{\tau \beta}{2} \mathbf{u}^T \left(D\otimes \mathcal{M}\right) \mathbf{u},
\end{equation}
where $\mathcal{M}$ and $\mathcal{M}_\gamma$ are defined in~(\ref{eq:def-Mcal}), $D$ in~(\ref{eq:CD}) below,
$\mathbf{y}$, $\mathbf{y}_d$, and $\mathbf{u}$ %and $\mathbf{g}$ 
are vectors corresponding to the state, desired state, and control, respectively,
%and boundary conditions respectively.
that contain concatenated vectors $\mathbf{y}_i, \mathbf{{y}_d}_i, \mathbf{u}_i \in \mathbb{R}^{n_h n_\xi \times 1}, i=1, \ldots, n_t$,
due to the time-stepping,
\[
    \mathbf{y}=
    \begin{bmatrix}
        \mathbf{y}_1 & \dots & \mathbf{y}_{n_t}
    \end{bmatrix}^T,\,
    \mathbf{y}_d=
    \begin{bmatrix}
        \mathbf{{y}_d}_1 & \dots & \mathbf{{y}_d}_{n_t}
    \end{bmatrix}^T,\,
    \text {and }
    \mathbf{u}=\begin{bmatrix}
        \mathbf{u}_1 & \dots & \mathbf{u}_{n_t}
    \end{bmatrix}^T.
\]
%we write 
%\begin{equation*}
% \mathbf{y}=\begin{bmatrix}
% \mathbf{y}_1 \\
% \vdots \\
% \mathbf{y}_{n_t}
% \end{bmatrix}, \quad 
% \mathbf{y}_d=\begin{bmatrix}
% \mathbf{{y}_d}_1 \\
% \vdots \\
% \mathbf{{y}_d}_{n_t}
% \end{bmatrix} , \quad 
% \text {and } \mathbf{u}=\begin{bmatrix}
% \mathbf{u}_1 \\
% \vdots \\
% \mathbf{u}_{n_t}
% \end{bmatrix} \text {. }
% \end{equation*}
%Considering the implicit Euler method to discretize in time time the PDE constraint 
%\begin{equation}
%    \label{eq: implicit euler scheme}
%    \frac{y^k-y^{k-1}}{\tau}-\Delta y^k=u^k,
%\end{equation}
%with $n_t$ time steps of size $\tau$ where $n_t=T/\tau$ represents the number of time steps over the interval $[0,T]$ and the finite element discretization of the weak form then gives 
After the application of the stochastic Galerkin finite element discretization to~(\ref{eq:constraint:stochastic-time-dependent}),
and using the implicit Euler method for the time discretization, we obtain
\begin{equation}
    \label{eq: back Euler}
    M\mathbf{y}_k+\tau A\mathbf{y}_k=M\mathbf{y}_{k-1}+\tau M \mathbf{u}_k.
\end{equation}
Combining all time steps of~(\ref{eq: back Euler}) in all-at-once discretization%of the state equation (\ref{eq:constraint:stochastic-time-dependent}) 
~(\cite{MR3023467,rees2010all}), we can write
% \begin{equation*}
%     \mathcal{A}_t \mathbf{y}-\tau \mathcal{N} \mathbf{u}=
%     %   \begin{bmatrix}
%     %    \mathcal{M} \mathrm{y}_0 &0& \cdots &0
%     %    \end{bmatrix}^T,
%     \begin{bmatrix}
%         \mathcal{M} \mathrm{y}_0 \, , & 0 \, , & \dots ,0
%     \end{bmatrix}^T.
%     %\mathrm{d},
% \end{equation*}
\zh{%
\begin{equation*}
    \mathcal{A}_t \mathbf{y}-\tau \mathcal{N} \mathbf{u}=
\begin{bmatrix}
        \mathcal{M} \mathrm{y}_0 + \mathbf{g} \, , & \mathbf{g} \, , & \dots ,\mathbf{g}
    \end{bmatrix}^T \eqqcolon \mathbf{d},
\end{equation*}
the vector $\mathbf{g}$ contains the contributions from the Dirichlet boundary data at each time step, and
}

\begin{equation}
    \label{eq: K_t for stochastic}
    \mathcal{A}_t=\left(I_{n_t}\otimes \mathcal{L}\right)-\left(C\otimes \mathcal{M}\right), \,
    \mathcal{N}=I_{n_t}\otimes H_1\otimes M,\text{ and }
    \mathcal{L}=H_1\otimes \left(M+\tau A_1\right)+\tau\sum_{\ell=2}^{n_A}H_\ell\otimes A_\ell.
\end{equation}
% \begin{equation}\label{A_t&N}
%     \mathcal{A}_t=\begin{bmatrix}
%     \mathcal{L} & & & \\
%     -\mathcal{M} & \mathcal{L} & & \\
%     & \ddots & \ddots & \\
%     & & -\mathcal{M} & \mathcal{L}
%     \end{bmatrix}, \, \qquad \,
%     \mathcal{N}=\begin{bmatrix}
%     \mathcal{M} & & & \\
%     & \mathcal{M} & & \\
%     & & \ddots & \\
%     & & & \mathcal{M}
%     \end{bmatrix}, \,
%     % \mathbf{d}=\begin{bmatrix}
%     % \mathcal{M} \mathrm{y}_0 \\
%     % 0 \\
%     % \vdots \\
%     % 0
%     % \end{bmatrix},
%     \end{equation}
% and
% \begin{equation}
%     \label{eq:Lcal}
%     \mathcal{L}=H_1\otimes \left(M+\tau A_1\right)+\tau\sum_{\ell=2}^{n_A}H_\ell\otimes A_\ell.
% \end{equation}
%Also, $\mathcal{N}$ can be rewritten as
%\begin{equation}
%    \label{eq:Ncal}
%    \mathcal{N}=I_{n_t}\otimes H_1\otimes M.
%\end{equation}
% The matrices $\mathcal{A}_t$ and $\mathcal{N}$ can be constructed using Kronecker product as
% \begin{equation}
%     \label{eq: K_t for stochastic}
%     \mathcal{A}_t=\left(I_{n_t}\otimes \mathcal{L}\right)-\left(C\otimes \mathcal{M}\right), \qquad  
%     \mathcal{N}=I_{n_t}\otimes H_1\otimes M, 
% \end{equation}
where the matrix~$C$ and matrix~$D$, used in~(\ref{eq:obj:stoch-time-dependent-discretize}) and also below, are defined as
\begin{equation}
    % C=\begin{bmatrix}
    % 0 & & & \\
    % -1 & 0 & & \\
    % & \ddots & \ddots & \\
    % & & -1 & 0
    % \end{bmatrix},
    % \quad
    % D=\begin{bmatrix}
    %     \frac12&&&&\\
    %     &1&&&\\
    %     &&\ddots&\\
    %     &&&1&\\
    %     &&&&\frac12
    % \end{bmatrix}.
    C=\operatorname{subdiag}(1, 1, \dots, 1)_{-1}, \quad
    D=\operatorname{diag}\left(\frac{1}{2}, 1, \dots, 1, \frac{1}{2}\right).
    \label{eq:CD}
\end{equation}
Forming the Lagrangean and applying the first-order optimality conditions, we get %the KKT system
\begin{equation}\label{eq:KKT:stoch-time-dependent}
    \begin{bmatrix}
        \tau D\otimes\mathcal{M}_{\gamma} & 0                              & -\mathcal{A}_t^T   \\
        0                                 & \beta \tau D\otimes\mathcal{M} & \tau \mathcal{N}^T \\
        -\mathcal{A}_t                    & \tau \mathcal{N}               & 0
    \end{bmatrix}
    \begin{bmatrix}
        \mathbf{y} \\
        \mathbf{u} \\
        \zhl{\boldsymbol{\lambda}}\end{bmatrix}
    =\begin{bmatrix}
        \tau \left(D\otimes\mathcal{M}\right)\cdot\left(\mathbf{1}_{n_t}\otimes \mathbf{y}_d\right) \\
        \mathbf{0}                                                                                  \\
        \mathbf{d}
    \end{bmatrix},
\end{equation}
where the boundary-inclusive source vector $\mathbf{d}$ is defined as above, and $\mathbf{1}_{n_t}\in \mathbb{R}^{n_t\times 1}$ is the column all-ones vector.
% where $\mathbf{d}=\begin{bmatrix}
%         \mathcal{M} \mathrm{y}_0+\mathbf{g} & \mathbf{g} & \cdots & \mathbf{g}
%     \end{bmatrix}^T$, and $\mathbf{1}_{n_t}\in \mathbb{R}^{n_t\times 1}$ is the column all-ones vector.

\subsection{PINT-based block-diagonal hierarchical Gauss-Seidel preconditioner}
In analogy to (\ref{eq:stoch-prec}), we propose a preconditioner for~(\ref{eq:KKT:stoch-time-dependent}) as %the inverse of 
\begin{equation}
    \label{eq: precon for time}
    %    \bar{\mathcal{P}}=
    \mathcal{P}_\text{hGSoc-PINT} \approx \begin{bmatrix}
        \tau D\otimes \mathcal{M}_\gamma &                                 &                   \\
                                         & \tau \beta D\otimes \mathcal{M} &                   \\
                                         &                                 & \bar{\mathcal{S}}
    \end{bmatrix},
\end{equation}
where $\bar{\mathcal{S}}$ is a computationally efficient approximation of the exact Schur complement
\begin{equation}
    \label{time-exact-pre}
    \bar{\mathcal{S}}_{\text{exact}}=\frac{1}{\tau}\mathcal{A}_t\left(D\otimes \mathcal{M}_\gamma\right)^{-1}
    \mathcal{A}_t^T+\frac{\beta}{\tau} \mathcal{N} \left(D\otimes \mathcal{M}\right)^{-1}\mathcal{N}^T.
\end{equation}
%The (1,1)-block is based Kronecker products of $D$, $H_1$, and $M$. The (2,2)-block is generated by Kronecker products of $\beta D$, $H_1$, and $M$. The (3,3)-block is a Schur complement of the (1,1) and (2,2) blocks.
%In practice, we can process these three blocks efficiently via parallel computing. 
% For example, for the (1,1)-block $\mathfrak{A}_1$, we can use the Kronecker product to process it efficiently. 
%
%The preconditioner is block diagonal, so the solves with the diagonal blocks can be performed independently. 
%
The first two blocks that correspond to (scaling of) the mass matrix are approximated
by Chebyshev semi-iteration. %from Algorithm~\ref{alg:chebmass}. 
Since the iteration entails matrix-vector multiplications, we note that %these can be written, 
using~(\ref{eq:mat-vec}) we have %as 
\[
\begin{aligned}
     \left(\tau D\otimes H_1\otimes M\right)\mathbf{v_1}&=
    %    \tau
    %    \begin{bmatrix}
    %        \frac12 H_\gamma\otimes M & & & \\
    %        &  H_\gamma\otimes M & & \\
    %        & & \ddots & \\
    %        & & & \frac12 H_\gamma\otimes M
    %    \end{bmatrix}
    %    \begin{bmatrix}
    %        \mathbf{v_{11}}\\
    %        \mathbf{v_{12}}\\
    %        \vdots\\
    %        \mathbf{v_{1n_t}}
    %    \end{bmatrix}=
    \tau
    \begin{bmatrix}
        \frac12 \vectorize{\left(MV_{11}H_\gamma\right)} &
        \vectorize{\left(MV_{12}H_\gamma\right)}         &
        \cdots                                           &
        \frac12 \vectorize{\left(MV_{1n_t}H_\gamma\right)}
    \end{bmatrix}^T           \\
    % \begin{bmatrix}
    %     \frac12 \vectorize{\left(MV_{11}H_\gamma\right)}\\
    %     \vectorize{\left(MV_{12}H_\gamma\right)}\\
    %     \vdots \\
    %     \frac12 \vectorize{\left(MV_{1n_t}H_\gamma\right)}\\
    % \end{bmatrix},\\
     \left(\beta\tau D\otimes H_1\otimes M\right)\mathbf{v_2}&=
    %    \beta\tau
    %    \begin{bmatrix}
    %        \frac12 H_1\otimes M & & & \\
    %        &  H_1\otimes M & & \\
    %        & & \ddots & \\
    %        & & & \frac12 H_1\otimes M
    %    \end{bmatrix}
    %    \begin{bmatrix}
    %        \mathbf{v_{21}}\\
    %        \mathbf{v_{22}}\\
    %        \vdots\\
    %        \mathbf{v_{2n_t}}
    %    \end{bmatrix}=
    \beta\tau
    \begin{bmatrix}
        \frac12 \vectorize{\left(MV_{21}H_1\right)} &
        \vectorize{\left(MV_{22}H_1\right)}         &
        \cdots                                      &
        \frac12 \vectorize{\left(MV_{2n_t}H_1\right)}
    \end{bmatrix}^T,
    % \begin{bmatrix}
    %     \frac12 \vectorize{\left(MV_{21}H_1\right)}\\
    %     \vectorize{\left(MV_{22}H_1\right)}\\
    %     \vdots \\
    %     \frac12 \vectorize{\left(MV_{2n_t}H_1\right)}\\
    % \end{bmatrix},
\end{aligned}\]
where $\mathbf{v_1}$ and $\mathbf{v_2}$ are the vectors %of size $n_t n_\xi n_h$ 
obtained by concatenating $\mathbf{v_{1i}}$ and $\mathbf{v_{2i}}$, $i=1,\dots,n_t$, respectively,
which correspond to the time steps.
The matrices~$V_{1i},\, V_{2i}\in \mathbb{R}^{n_h\times n_\xi}$
are then the \emph{matricized} counterparts of~$\mathbf{v_{1i}}$ and $\mathbf{v_{2i}}$, respectively.
Next, since inverting~$\bar{\mathcal{S}}_{\text{exact}}$ is computationally prohibitive, we propose an approximation as
\begin{equation}
    \label{eq: S hat}
    \bar{\mathcal{S}}=\frac{1}{\tau}\underbrace{\left(\mathcal{A}_t+\tau\sqrt{\frac{1+\gamma}{\beta}}\mathcal{N}\right)}_{\eqqcolon \bar{\mathcal{Z}}}\left(D\otimes \mathcal{M}_\gamma\right)^{-1}
    \left(\mathcal{A}_t+\tau\sqrt{\frac{1+\gamma}{\beta}}\mathcal{N}\right)^T,
\end{equation}
where using \eqref{eq: K_t for stochastic}, we can rewrite $\bar{\mathcal{Z}}$ as
\begin{equation}
    \label{Z hat}
    \bar{\mathcal{Z}}=I_{n_t}\otimes\left(H_1\otimes\left( M+\tau A_1\right)+\tau\sum_{\ell=2}^{n_A} H_\ell\otimes A_\ell\right)+\tau\sqrt{\frac{1+\gamma}{\beta}}I_{n_t}\otimes H_1\otimes M-C\otimes H_1\otimes M.
    % &=\left(I_{n_t}+C\right)\otimes\left(H_1\otimes M\right)+\tau I_{n_t}\otimes\left(\sum_{i=0}^N H_\ell\otimes A_\ell\right)+\frac{\tau}{\sqrt{\beta}}\left(I_{n_t}\otimes H_1 \otimes M\right)\\
    %&=\sum_{i=0}^N\underbrace{\left(I_{n_t}\otimes G_l\right)}_{\Gcal_i}\otimes \hat{A}_i+\underbrace{C\otimes G_0}_{\hat{C}}\otimes M,\\
    %&=\sum_{i=0}^N\underbrace{\left(I_{n_t}\otimes G_l\right)}_{\Gcal_i}\otimes \hat{A}_i+\underbrace{C\otimes G_0}_{\hat{C}}\otimes M,
\end{equation}
By dropping the last term~$C\otimes H_1\otimes M$, we further approximate $\bar{\mathcal{Z}}$ by %$\tilde{\mathcal{Z}}$, where
\begin{equation}
    \label{eq: Z hat approx}
    \tilde{\mathcal{Z}}= I_{n_t}\otimes\left\{ H_1\otimes \left[(1+\tau\sqrt{\frac{1+\gamma}{\beta}})M+\tau A_1\right]+\tau \sum_{\ell=2}^{n_A}H_\ell\otimes A_\ell\right\}.
\end{equation}
Observe that $\tilde{\mathcal{Z}}$ is symmetric, and it can be rewritten as 
\begin{equation}
\label{another form of Ztilde}
\tilde{\mathcal{Z}}=\tau \sqrt{\frac{1+\gamma}{\beta}}\mathcal{N} + I_{n_t}\otimes \mathcal{L}.
\end{equation}
The idea is to~use%
\zhl{\begin{equation}
        \label{S_tilde}
        \tilde{\mathcal{S}} \approx \tilde{\mathcal{Z}} \left( D \otimes \mathcal{M}_\gamma \right)^{-1} \tilde{\mathcal{Z}},
    \end{equation}}%
and in particular the solves with $\tilde{\mathcal{Z}}$ are approximated by %employing %
%the hierarchical Gauss-Seidel preconditioner from 
Algorithm~\ref{alg:hGS}, similarly to the steady-state~case.
We remark that by dropping all terms with $\ell>1$ from~$\tilde{\mathcal{Z}}$ in~(\ref{eq: Z hat approx}); that is, considering
\[
    \tilde{\mathcal{Z}}_0= I_{n_t}\otimes H_1\otimes \left[(1+\tau\sqrt{\frac{1+\gamma}{\beta}})M+\tau A_1\right],
\]
we recover the mean-based preconditioner~\cite{Benner-2016-BDP}. Since the application of Algorithm~\ref{alg:hGS}
entails matrix-vector multiplications, using~(\ref{eq:mat-vec}) we formulate $\tilde{\mathcal{Z}}\mathbf{v_3}$ as
\begin{equation}
    \begin{aligned}
        \label{eq: Zhat v_3}
        \tilde{\mathcal{Z}}\mathbf{v_3} & =\left(I_{n_t}\otimes\left\{H_1\otimes \left[(1+\tau\sqrt{\frac{1+\gamma}{\beta}})M+\tau A_1\right]+\tau \sum_{\ell=2}^{n_A} H_\ell\otimes A_\ell\right\}\right)\mathbf{v_3} \\
        %    &=\begin{bmatrix}
        %        \sum_{\ell=1}^{n_A} H_\ell\otimes \hat{A}_\ell&&&\\
        %        &\sum_{\ell=1}^{n_A} H_\ell\otimes \hat{A}_\ell&&\\
        %        &&\ddots&\\
        %        &&&\sum_{\ell=1}^{n_A} H_\ell\otimes \hat{A}_\ell
        %    \end{bmatrix}
        %    \begin{bmatrix}
        %        \mathbf{v_{31}}\\
        %        \mathbf{v_{32}}\\
        %        \vdots\\
        %        \mathbf{v_{3n_t}} 
        %    \end{bmatrix}
                                        & =
        % \begin{bmatrix}
        %     \sum_{\ell=1}^{n_A} \vectorize{\left(\hat{A}_\ell V_{31}H_\ell\right)}\\
        %     \sum_{\ell=1}^{n_A} \vectorize{\left(\hat{A}_\ell V_{32}H_\ell\right)}\\
        %     \vdots\\
        %     \sum_{\ell=1}^{n_A} \vectorize{\left(\hat{A}_\ell V_{3n_t}H_\ell\right)}
        % \end{bmatrix}
        \begin{bmatrix}
            \sum_{\ell=1}^{n_A} \vectorize{\left(\hat{A}_\ell V_{31}H_\ell\right)} &
            \sum_{\ell=1}^{n_A} \vectorize{\left(\hat{A}_\ell V_{32}H_\ell\right)} &
            \cdots                                                                 &
            \sum_{\ell=1}^{n_A} \vectorize{\left(\hat{A}_\ell V_{3n_t}H_\ell\right)}
        \end{bmatrix}^T
        ,
    \end{aligned}
\end{equation}
where
\begin{equation}
    \label{eq:A_hat}
    \hat{A}_\ell=
    \begin{cases}
        \left(1+\tau\sqrt{\frac{1+\gamma}{\beta}}\right)M+\tau A_1, & \ell=1,           \\
        \tau A_\ell,                                                & \ell=2,\dots,n_A,
    \end{cases}
\end{equation}
and $V_{3i}\in \mathbb{R}^{n_h\times n_\xi}$ is the \emph{matricized} form of
the $i$-th block of $\mathbf{v_3}$, with $i=1,\dots,n_t$.
%$\mathbf{v_{3i}}$ and $\mathbf{v_3} \in \mathbb{R}^{n_t n_\xi}$

The practical implementation of the preconditioner for the time-dependent system leverages the inherent structure of the
\emph{all-at-once} formulation.
As defined in \eqref{eq: Z hat approx}, the core operator of the Schur complement preconditioner, $\tilde{\mathcal{Z}}$,
is block-diagonal with respect to the time steps. This structure extends to the entire KKT system, which then
makes the preconditioning easily parallelizable.
Specifically, an application of the preconditioner $\mathcal{P}_\text{hGSoc-PINT}$ entails
an application of the steady-state optimal control
preconditioner $\mathcal{P}_\text{hGSoc}$ from Algorithm~\ref{alg:hGSoc-1}--\ref{alg:hGSoc-2} to all time steps simultaneously,
and so it represents \emph{parallel-in-time} (PINT) approach. It is summarized as Algorithm~\ref{alg:hGSoc-PINT}.

\begin{algorithm}[H]
    \caption{Parallel-in-time hGSoc preconditioner (hGSoc-PINT)}
    The preconditioner  $\mathcal{P}_\text{hGSoc-PINT}: \bar{\mathbf{R}} \longmapsto \bar{\mathbf{V}}$ is defined as:
    \begin{algorithmic}[1]
        \For{$k = 1, \dots, n_t$}
        \State Extract $(\bar{R}^Y_k, \bar{R}^U_k, \bar{R}^\Lambda_k)$. \hfill (the subvector~$k$ of $\bar{\mathbf{R}}$) % for the time step $k$)
        \State Calculate $(\bar{V}^Y_k, \bar{V}^U_k, \bar{V}^\Lambda_k) = \mathcal{P}_{hGSoc}(\bar{R}^Y_k, \bar{R}^U_k, \bar{R}^\Lambda_k)$ \hfill (apply Algorithm~\ref{alg:hGSoc-1}--\ref{alg:hGSoc-2})
        \EndFor
        \State Concatenate $\{(\bar{V}^Y_k, \bar{V}^U_k, \bar{V}^\Lambda_k)\}_{k=1}^{n_t}$ into $\bar{\mathbf{V}}$.
    \end{algorithmic}
    \label{alg:hGSoc-PINT}
\end{algorithm}

\section{Spectral analysis of the preconditioners}
\label{sec: spectral analysis}
Since the steady-state optimal control problem can be viewed as a special case of the time-dependent formulation (with $n_t=1$), we focus on analyzing the time-dependent setting; the steady-state results then follow as a direct consequence.
The all-at-once discretization, presented in Section~\ref{sec:stochastic-model-td},
couples all time steps simultaneously, yielding a significantly larger KKT system than its steady-state counterpart.
Our goal is to prove that the proposed parallel-in-time preconditioner,
based on the hGSoc-PINT in Algorithm~\ref{alg:hGSoc-PINT},
is spectrally equivalent to the ideal (but computationally more expensive) preconditioner.

\begin{definition}[Spectral Equivalence]
    Two symmetric positive definite (SPD) matrices $A$ and $B$ are said to be spectrally equivalent, denoted $A\sim B$, if there exist positive constants $a\leq b$, such that
    $$
        a \mathbf{v}^T B \mathbf{v} \le \mathbf{v}^T A \mathbf{v} \le b \mathbf{v}^T B \mathbf{v}
    $$
    holds for all non-zero vectors $\mathbf{v}$. Equivalently, all eigenvalues of the preconditioned matrix $B^{-1}A$ are contained within the fixed interval, which means $\lambda(B^{-1}A)\subset [a, b]$.
\end{definition}
The proof proceeds by establishing spectral equivalences
\begin{equation}
    \label{eq:spectral-equiv}
    \bar{\mathcal{S}}_{\text{exact}}\sim \bar{\mathcal{S}} \sim \tilde{\mathcal{S}} \sim \tilde{\mathcal{S}}_r\sim\tilde{\mathcal{S}}_\text{hGS-PINT}.
\end{equation}
% Here, $\bar{\mathcal{S}}_{\text{exact}}$ denotes the exact Schur complement~\eqref{time-exact-pre},
% $\bar{\mathcal{S}}$ \zhl{defined in~\eqref{eq: S hat}}, is an approximation \zhl{of $\bar{\mathcal{S}}_{\text{exact}}$}, 
% and $\tilde{\mathcal{S}}$ \zhl{defined in \eqref{S_tilde}} is a block-diagonal approximation \zhl{of $\bar{\mathcal{S}}$}.
\zhl{
    Here, $\bar{\mathcal{S}}_{\text{exact}}$ denotes the exact Schur complement~\eqref{time-exact-pre}. The matrix $\bar{\mathcal{S}}$, defined in~\eqref{eq: S hat}, serves as an approximation of $\bar{\mathcal{S}}_{\text{exact}}$, while $\tilde{\mathcal{S}}$ given in~\eqref{S_tilde} is a block-diagonal approximation of $\bar{\mathcal{S}}$.}
% obtained by replacing $\bar{\mathcal{Z}}$ in~\eqref{eq: S hat} with $\tilde{\mathcal{Z}}$ from~\eqref{eq: Z hat approx}, which eliminates time-coupling terms and enables parallel-in-time computation.
The operator $\tilde{\mathcal{S}}_r$ represents the truncated hierarchical operator
    \begin{equation}
    \label{S_r}    
    \tilde{\mathcal{S}}_r \;=\; \tilde{\mathcal{Z}}_r\,(D\otimes \mathcal{M}_\gamma)^{-1}\,\tilde{\mathcal{Z}}_r^{T},
    \,
    \tilde{\mathcal{Z}}_r \;=\; I_{n_t}\otimes\left\{ H_1\otimes \left[(1+\tau\sqrt{\tfrac{1+\gamma}{\beta}})M+\tau A_1\right]+\tau \sum_{\ell=2}^{r}H_\ell\otimes A_\ell\right\},
    \end{equation}
            with $r=1,\dots,n_A$. When $r=1$, $\tilde{\mathcal{S}}_r$ reduces to the mean-based preconditioner employed in~\cite{Benner-2016-BDP}, and when $r=n_A$, it recovers the full operator $\tilde{\mathcal{S}}$.
Finally, $\tilde{\mathcal{S}}_{\text{hGS-PINT}}=\tilde{\mathcal{Z}}_{\text{hGS-PINT}}(D\otimes\mathcal{M}_\gamma)^{-1}\tilde{\mathcal{Z}}_{\text{hGS-PINT}}^T$ represents the computationally feasible approximation of $\tilde{\mathcal{S}}_r$ in which the linear systems $\tilde{\mathcal{Z}}_r \mathbf{x} = \mathbf{b}$ are solved approximately via the hierarchical Gauss-Seidel method (Algorithm~\ref{alg:hGS}), as implemented in the parallel-in-time framework of Algorithm~\ref{alg:hGSoc-PINT}.

To establish the spectral equivalences in this chain, we require several technical results. We state some auxiliary lemmas concerning matrix perturbations and congruence transformations, which will serve as building blocks for the main theorems. Their proofs are provided in Appendix~\ref{sec:appendix}. In what follows, we denote by $\sigma_{\min}(\cdot)$ and $\sigma_{\max}(\cdot)$ the smallest and largest singular values, respectively.%

\begin{lemma}[Spectral Equivalence via Error Bound]
    \label{lem:error_bound}
    Let $A$ and $B$ be symmetric positive definite matrices. Define the error matrix $E = A - B$. If there exists a constant $0<\delta < 1$ such that for all non-zero vectors $\mathbf{v}$,
    $$|\mathbf{v}^T E \mathbf{v}| \le \delta (\mathbf{v}^T B \mathbf{v}),$$
    then~$A$ and~$B$ are spectrally equivalent with 
    $$1 - \delta \le \lambda(B^{-1}A) \le 1 + \delta.$$
\end{lemma}
% \begin{proof}
% See Appendix~\ref{sec:appendix} for the proof.
% \end{proof}

\begin{lemma}[Eigenvalues under Congruence Transformation]
    \label{lem:congruence}
    Let $C$ and $D$ be symmetric positive definite matrices, and let $Q$ be a nonsingular matrix. The eigenvalues of the pair $(C, D)$ are identical to those of the transformed pair $(Q^T C Q, Q^T D Q)$.
\end{lemma}
% \begin{proof}
%     See Appendix~\ref{sec:appendix} for the proof.
% \end{proof}
\zhl{
    \begin{lemma}
        \label{lem:sig<=l2}
        For any nonzero vector $\mathbf{v}$, we have
        $            \sigma_{\min}(A)\|\mathbf{v}\|_2 \leq \|A\mathbf{v}\|_2.$
    \end{lemma}}%
% \begin{proof}
%     See Appendix~\ref{sec:appendix} for the proof.
% \end{proof}
\begin{lemma}\label{lemma ineq}
    For any matrices $A$ and $B$ of the same dimensions, %the following inequality holds:
    \[
        \sigma_{\min}(A+B)\geq \sigma_{\min}(B)-\|A\|_2.
    \]
    %where $\sigma_{\min}(\cdot)$ is the smallest singular value of its matrix argument.
\end{lemma}

With these auxiliary results in place, we now proceed to establish the spectral equivalences in~\eqref{eq:spectral-equiv}. We begin by recalling the following spectral equivalence result from~\cite{Benner-2016-BDP}.
%, which shows that our factorized approximation~\eqref{eq: S hat} is spectrally close to the exact Schur complement.
\begin{theorem}[\cite{Benner-2016-BDP}, Theorems 4, 6] \label{thm:approx}
    Let $\bar{\mathcal{S}}_{\text{exact}}$ be the exact Schur complement and $\bar{\mathcal{S}}$ be its approximation as given by~\eqref{time-exact-pre} and~\eqref{eq: S hat}, respectively.    %For any $\alpha \ge 0$ satisfying the condition below, 
    The eigenvalues of $\bar{\mathcal{S}}^{-1}\bar{\mathcal{S}}_{\text{exact}}$ are given by
\[\lambda(\bar{\mathcal{S}}^{-1}\bar{\mathcal{S}}_{\text{exact}}) \subseteq \left[ \frac{1}{2(1+\alpha)}, 1 \right),
\]
where $\alpha$ satisfies
    $        \alpha < \left(\frac{\sqrt{\kappa(\mathcal{A}_t)}+1}{\sqrt{\kappa(\mathcal{A}_t)}-1}\right)^2-1,
$
    and $\mathcal{A}_t$ is defined in~\eqref{eq: K_t for stochastic}.
    % depends on the problem type:
    % \begin{itemize}
    %     \item For the steady-state problem, where $\hat{\mathcal{S}} = \mathcal{S}$ from~\eqref{eq:stoch-prec}:
    %     \[ \alpha < \left(\frac{\sqrt{\kappa(\mathcal{M}^{-\frac{1}{2}}\mathcal{A}\mathcal{M}^{-\frac{1}{2}})}+1}{\sqrt{\kappa(\mathcal{M}^{-\frac{1}{2}}\mathcal{A}\mathcal{M}^{-\frac{1}{2}})}-1}\right)^2-1. \]
    %     \item For the time-dependent problem, where $\hat{\mathcal{S}} = \overline{\mathcal{S}}$ from~\eqref{eq: precon for time}:
    %     \[ \alpha < \left(\frac{\sqrt{\kappa(\mathcal{A}_t)}+1}{\sqrt{\kappa(\mathcal{A}_t)}-1}\right)^2-1. \]
    % \end{itemize}
\end{theorem}
Next, we show the relationship between $\bar{\mathcal{S}}$ and~$\tilde{\mathcal{S}}$, as given by \eqref{eq: S hat} and~\eqref{S_tilde}, respectively. First, however, we prove the following lemma.
\begin{lemma}
    \label{lem:lower_bound}
    % Under certain assumptions, the minimum eigenvalue of $\tilde{\mathcal{S}}$ has the following lower bound 
    % {\color{blue}(bs: What assumptions?)}
    % {\color{blue}(zl: I rewrote this part in blue}
    Let $W=(D\otimes \mathcal{M}_\gamma)^{-1}$, and assume there exists a constant $\mu > 1$ such that
    \begin{equation}
        \label{mucond}
        \tau \sqrt{\frac{1+\gamma}{\beta}}\sigma_{\min} (\mathcal{N} W^\frac12) \geq \mu \|(I_{n_t}\otimes \mathcal{L})W^\frac12\|_2,
    \end{equation}where $\mathcal{L}$ is defined in~\eqref{eq: K_t for stochastic}.
    Then the minimum eigenvalue of $\tilde{\mathcal{S}}$ has the following lower bound
    \[
        \lambda_{\min}(\tilde{\mathcal{S}}) \geq \frac{\tau}{\beta}\left(1-\frac{1}{\mu}\right)^2 \sigma^2_{\min}(M^{\frac12}).
    \]
\end{lemma}
\begin{proof}
    Using the definition of $\tilde{\mathcal{S}}$ in~\eqref{S_tilde}, and the fact that for any real matrix~$A$, $\lambda(AA^T)=\sigma^2_{\min}(A)$, we have 
    $$\lambda_{\min}(\tilde{\mathcal{S}}) = \lambda_{\min}(\tilde{\mathcal{Z}}W\tilde{\mathcal{Z}}^T) =\lambda_{\min}((\tilde{\mathcal{Z}}W^{\frac12})(\tilde{\mathcal{Z}}W^{\frac12})^T)=
    \sigma^2_{\min}(\tilde{\mathcal{Z}}W^{\frac12}).
$$
    Recall from~\eqref{another form of Ztilde} that $\tilde{\mathcal{Z}}W^\frac12 = \tau \sqrt{\frac{1+\gamma}{\beta}}\mathcal{N} W^{\frac12} + (I_{n_t}\otimes \mathcal{L})W^\frac12$. Now, applying Lemma~\ref{lemma ineq} yields
    \begin{equation}
    \label{4.3}
        \sigma_{\min}(\tilde{\mathcal{Z}}W^\frac12) \geq \tau \sqrt{\frac{1+\gamma}{\beta}}\sigma_{\min} (\mathcal{N} W^\frac12) - \|(I_{n_t}\otimes \mathcal{L})W^\frac12\|_2.
    \end{equation}
    Now using~\eqref{4.3} and the assumption~\eqref{mucond}, we obtain 
    \footnote{In our experience, this condition is often satisfied numerically when $\tau \gg \sqrt{\beta}$.}
    % We assume there exists a constant \mu > 1 such that ...
    % Using the assumption (\ref{mucond}), we obtain
    % Since condition (\ref{mucond}) holds by hypothesis,
    \[
        \sigma_{\min}(\tilde{\mathcal{Z}}W^\frac12) \geq \left(1-\frac{1}{\mu}\right) \tau \sqrt{\frac{1+\gamma}{\beta}} \sigma_{\min}(\mathcal{N} W^\frac12).
    \]
    Consequently, the lower bound for $\lambda_{\min}(\tilde{\mathcal{S}})$ is
    \[
        \lambda_{\min}(\tilde{\mathcal{S}}) \geq \left(1-\frac{1}{\mu}\right)^2 \left( \tau \sqrt{\frac{1+\gamma}{\beta}}\right)^2 \sigma^2_{\min}(\mathcal{N} W^{\frac12}).
    \]
    Next, observe that
    \[
        \sigma^2_{\min}(\mathcal{N} W^{\frac12}) = \sigma^2_{\min}\left( (I_{n_t}\otimes\mathcal{M}) (D^{-\frac12}\otimes \mathcal{M}_\gamma^{-\frac12}) \right)
        = \sigma^2_{\min}(D^{-\frac12}) \cdot \sigma^2_{\min}(\mathcal{M} \mathcal{M}_\gamma^{-\frac12}),
    \]
    % \[
    %     \begin{aligned}
    %         \sigma^2_{\min}(\mathcal{N} W^{\frac12}) &= \sigma^2_{\min}\left( (I_{n_t}\otimes\mathcal{M}) (D^{-\frac12}\otimes \mathcal{M}_\gamma^{-\frac12}) \right) \\
    %         &= \sigma^2_{\min}\left( (I_{n_t} D^{-\frac12}) \otimes (\mathcal{M} \mathcal{M}_\gamma^{-\frac12}) \right) \\
    %         &= \sigma^2_{\min}(D^{-\frac12}) \cdot \sigma^2_{\min}(\mathcal{M} \mathcal{M}_\gamma^{-\frac12}),
    %     \end{aligned}
    %     \]
    where the last equality follows from the property that the singular values of a Kronecker product are the products of the singular values of the factors, i.e., $\sigma(A \otimes B) = \sigma(A)\sigma(B)$.

    It is easy to verify that
    $\sigma_{\min}(D^{-\frac12}) = 1$ and
    \[
        \sigma^2_{\min}(\mathcal{M} \mathcal{M}_\gamma^{-\frac12})
        =\sigma^2_{\min}( (H_1\otimes M)(H_\gamma^{-\frac{1}{2}}\otimes M^{-\frac{1}{2}}) )
        = \frac{1}{1+\gamma}\sigma^2_{\min}(M^\frac12).
    \]
    % \[
    % \begin{aligned}
    % \sigma^2_{\min}(\mathcal{M} \mathcal{M}_\gamma^{-\frac12}) 
    % &=\sigma^2_{\min}( (H_1\otimes M)(H_\gamma^{-\frac{1}{2}}\otimes M^{-\frac{1}{2}}) )\\
    % &=\sigma^2_{\min}( (H_1H_\gamma^{-\frac{1}{2}})\otimes(MM^{-\frac{1}{2}}) )\\
    % &= \frac{1}{1+\gamma}\sigma^2_{\min}(M^\frac12).
    % \end{aligned}
    % \]
    Thus, we have
    \[
        \sigma^2_{\min}(\mathcal{N} W^{\frac12}) = \frac{1}{1+\gamma} \sigma^2_{\min}(M^\frac12).
    \]
    Substituting this into the expression for the lower bound of $\lambda_{\min}(\tilde{\mathcal{S}})$, we finally get
    % \begin{equation} \label{lemma 2}
    \[
        \lambda_{\min}(\tilde{\mathcal{S}}) \geq \left(1-\frac{1}{\mu}\right)^2 \frac{\tau^2(1+\gamma)}{\beta} \left(\frac{1}{1+\gamma} \sigma^2_{\min}(M^\frac12)\right) = \frac{\tau}{\beta}\left(1-\frac{1}{\mu}\right)^2 \sigma^2_{\min}(M^{\frac12}).
    \]
    % \end{equation}
\end{proof}
%We can now state the following result.
\begin{theorem}
    Assume that   the conditions of
    Lemma~\ref{lem:lower_bound} hold.
    % {\color{blue}(zl: Assume $\tau$ is sufficiently large (e.g.,$\tau\gg \sqrt{\beta}\kappa(M^{\frac{1}{2}})$))}
    % {\color{blue}{zl: what about using $\tau\gg \sqrt{\beta}$ since $\kappa(M^{\frac{1}{2}}$ can be bounded with proper basis selection}} of %the preceding lemma 
    % Lemma~\ref{lem:lower_bound} hold. 
    Then, %spectral relationship between the matrices $\bar{\mathcal{S}}$ and $\tilde{\mathcal{S}}$ is bounded by 
    the eigenvalues $\lambda(\bar{\mathcal{S}}^{-1}\tilde{\mathcal{S}})$ satisfy
    \[
        (1-\theta)^2 \leq \lambda(\bar{\mathcal{S}}^{-1}\tilde{\mathcal{S}}) \leq (1+\theta)^2,
    \]
    where $\theta$  is the perturbation parameter, with
    \[
        \theta \coloneqq \sup_{\mathbf{v}\neq 0} \frac{\|W^\frac12(C\otimes\mathcal{M})^T\mathbf{v}\|_2}{\|W^\frac12\tilde{\mathcal{Z}}^T\mathbf{v}\|_2}\leq \frac{\sqrt{2\beta}}{\sqrt{1+\gamma}\left(1-\frac{1}{\mu}\right)\tau} \kappa(M^\frac12).
    \]
    % That is, for a proper parameter setting, $\bar{\mathcal{S}}$ is a good spectral approximation of $\tilde{\mathcal{S}}$.
\end{theorem}
\begin{proof}
    \zhl{
        Recall from~\eqref{Z hat} and~\eqref{eq: Z hat approx}, $\bar{\mathcal{Z}}=\tilde{\mathcal{Z}}-C\otimes\mathcal{M}$, which implies for any $\mathbf{v}\neq 0$,
        \[
            W^\frac12\bar{\mathcal{Z}}^T\mathbf{v} = W^\frac12\tilde{\mathcal{Z}}^T\mathbf{v} - W^\frac12(C\otimes\mathcal{M})^T\mathbf{v}.
        \]

        Using the triangle inequality, we have
        \[
            \left|\|W^\frac12\tilde{\mathcal{Z}}^T\mathbf{v}\|_2 - \|W^\frac12(C\otimes\mathcal{M})^T\mathbf{v}\|_2\right| \leq \|W^\frac12\bar{\mathcal{Z}}^T\mathbf{v}\|_2 \leq \|W^\frac12\tilde{\mathcal{Z}}^T\mathbf{v}\|_2 + \|W^\frac12(C\otimes\mathcal{M})^T\mathbf{v}\|_2.
        \]
        From the definition of $\theta$, we have
        $\|W^\frac12(C\otimes\mathcal{M})^T\mathbf{v}\|_2 \leq \theta \|W^\frac12\tilde{\mathcal{Z}}^T\mathbf{v}\|_2$. Substituting this into the inequality above, we obtain
        \[
            (1-\theta)\|W^\frac12\tilde{\mathcal{Z}}^T\mathbf{v}\|_2 \leq \|W^\frac12\bar{\mathcal{Z}}^T\mathbf{v}\|_2 \leq (1+\theta)\|W^\frac12\tilde{\mathcal{Z}}^T\mathbf{v}\|_2.
        \]
        Squaring these inequalities leads to the spectral bounds
        %we can bound $\lambda(\bar{\mathcal{S}}^{-1}\tilde{\mathcal{S}})$ as follows 
        \[
            (1-\theta)^2 \leq \lambda(\bar{\mathcal{S}}^{-1}\tilde{\mathcal{S}}) \leq (1+\theta)^2.
        \]
        Next, it remains to analyze the upper bound of $\theta$. To this end, using Lemma~\ref{lem:sig<=l2}, we know that
        $\sigma_{\min}(W^\frac12 \tilde{\mathcal{Z}}^T) \|\mathbf{v}\|_2\leq \|W^\frac12 \tilde{\mathcal{Z}}^T\mathbf{v}\|_2$.
        \[
            \theta = \sup_{\mathbf{v}\neq 0} \frac{\|W^\frac12(C\otimes\mathcal{M})^T\mathbf{v}\|_2}{\|W^\frac12\tilde{\mathcal{Z}}^T\mathbf{v}\|_2} \leq \frac{\|W^\frac12(C\otimes\mathcal{M})^T\|_2\|\mathbf{v}\|_2}{\sigma_{\min}(W^\frac12 \tilde{\mathcal{Z}}^T)\|\mathbf{v}\|_2}=\frac{\|W^\frac12(C\otimes\mathcal{M})^T\|_2}{\sigma_{\min}(W^\frac12 \tilde{\mathcal{Z}}^T)},
        \]
        % from which we bound the numerator and the denominator separately. First, for the numerator
        For the numerator, we have
        \[
            \begin{aligned}
                \|W^\frac12 (C\otimes \mathcal{M})^T\|_2 & = \|(D^{-\frac12}\otimes \mathcal{M}_\gamma^{-\frac12}) (C^T\otimes \mathcal{M}^T)\|_2 \\
                                                         & = \|(D^{-\frac12} C^T) \otimes (\mathcal{M}_\gamma^{-\frac12} \mathcal{M}^T)\|_2       \\
                                                         & = \|D^{-\frac12} C^T\|_2 \cdot \|\mathcal{M}_\gamma^{-\frac12} \mathcal{M}^T\|_2.
            \end{aligned}
        \]
        Furthermore, the spectral norm of the Kronecker product satisfies
        \[
            \|D^{-\frac12} C^T\|_2 \cdot \|\mathcal{M}_\gamma^{-\frac12} \mathcal{M}^T\|_2 \leq \frac{\sqrt{2}}{\sqrt{1+\gamma}}\|M^{\frac12}\|_2,
        \]
        which follows from the bound $\|D^{-1/2}C^T\|_2 \le 1$ and the fact that $\|\mathcal{M}_\gamma^{-1/2}\mathcal{M}^T\|_2$ is bounded by the scaled norm of the square root of the mass matrix.
        Thus, the numerator satisfies
        \[
            \|W^\frac12 (C\otimes \mathcal{M})^T\|_2 \leq \frac{\sqrt{2}}{\sqrt{1+\gamma}}\|M^{\frac12}\|_2.
        \]
        For the denominator, using lemma~\ref{lem:lower_bound}, we have
        \[
            \sigma_{\min}(W^\frac12 \tilde{\mathcal{Z}}^T) = \sigma_{\min}(\tilde{\mathcal{Z}} W^\frac12) \geq \left(1-\frac{1}{\mu}\right) \frac{\tau}{\sqrt{\beta}}\sigma_{\min}(M^\frac12).
        \]
        Combining the bounds for the numerator and denominator, %we obtain an upper bound for $\theta$ 
        \[
            \begin{aligned}
                \theta & \leq \frac{\frac{\sqrt{2}}{\sqrt{1+\gamma}}\|M^{\frac12}\|_2}{\left(1-\frac{1}{\mu}\right) \frac{\tau}{\sqrt{\beta}}\sigma_{\min}(M^\frac12)} \\
                       & = \frac{\sqrt{2\beta}}{\sqrt{1+\gamma}\left(1-\frac{1}{\mu}\right)\tau} \cdot \frac{\sigma_{\max}(M^\frac12)}{\sigma_{\min}(M^\frac12)}       \\
                       & = \frac{\sqrt{2\beta}}{\sqrt{1+\gamma}\left(1-\frac{1}{\mu}\right)\tau} \kappa(M^\frac12),
            \end{aligned}
        \]
        thereby completing the proof of the theorem.}
    % This upper bound indicates that when the parameter $\tau$ is sufficiently large (e.g., $\tau \gg \sqrt{\beta} \kappa(M^\frac12)$), $\theta$ will be a small quantity less than 1. This ensures that $\bar{\mathcal{S}}$ and $\tilde{\mathcal{S}}$ are spectrally close.
\end{proof}

Note that the  bound on $\theta$ shows that for sufficiently large~$\tau$ (e.g.,~$\tau \gtrsim \sqrt{\beta}\,\kappa(M)^{1/2}$), we have $\theta<1$; hence, $\bar{\mathcal S}$ and $\tilde{\mathcal S}$ are spectrally equivalent. In practice, $\kappa(M)$ can be kept $\mathcal{O}(1)$ via appropriate basis choices and mass lumping, so it suffices to require $\tau \gg \sqrt{\beta}$.

\begin{color}{black}
    Next, we establish the spectral equivalence between $\tilde{\mathcal{S}}$ and its truncated form $\tilde{\mathcal{S}}_r$. %Instead of relying on the spectral bounds of the factors, 
    To this end, we shall rely on the error matrix $\mathcal{E}_r = \tilde{\mathcal{S}} - \tilde{\mathcal{S}}_r$, %Based on the decay of the coefficients in the polynomial chaos expansion, % (assuming standard summability conditions on the coefficients $\gamma_\ell$), the error is dominated by the tail of the expansion.
    and the result of Lemma~\ref{lem:error_bound}.%, we need to bound the ratio $|\mathbf{v}^T \mathcal{E}_r \mathbf{v}| / |\mathbf{v}^T \tilde{\mathcal{S}}_r \mathbf{v}|$.

    \begin{theorem}[Truncation Error Bound]
        \label{thm:truncation}
        Let $\tilde{\mathcal{S}}$ be the full operator and $\tilde{\mathcal{S}}_r$ be the truncated hierarchical operator with truncation parameter $r$ defined in~\eqref{S_tilde} and~\eqref{S_r}, respectively.
        Assume $\nu \coloneqq \lambda_{\min}(\hat{A}_1) > 0$. %that the mean-field part of the operator $\tilde{\mathcal{Z}}$ is coercive, i.e., let$\nu \coloneqq \lambda_{\min}(\hat{A}_1) > 0$.
        %Assume that the mean-field part of the operator is coercive, i.e., there exists $\nu > 0$ such that $\mathbf{v}^T \tilde{\mathcal{S}}_1 \mathbf{v} \ge \nu \|\mathbf{v}\|^2$ for all $\mathbf{v}$.
        Define the truncation tail sum as
        \[
            \delta_r \coloneqq \frac{1}{\nu} \sum_{\ell=r+1}^{n_A} \|H_\ell\|_2 \, \|\hat{A}_\ell\|_2.
        \]
        If the decay of the stochastic expansion coefficients is such that $0<\delta_r < 1$, then~$\tilde{\mathcal{S}}$ and $\tilde{\mathcal{S}}_r$ are spectrally equivalent
        \[
            1-\delta_r \leq \lambda(\tilde{\mathcal{S}}_r^{-1}\tilde{\mathcal{S}}) \leq 1+\delta_r.
        \]
    \end{theorem}

    \begin{proof}
        Consider the error matrix $\mathcal{E}_r = \tilde{\mathcal{S}} - \tilde{\mathcal{S}}_r$. Based on the expansion structure of the operator defined in \eqref{eq: Z hat approx}, the error term consists of the neglected high-order terms in the stochastic expansion:
        \[
            \mathcal{E}_r = \sum_{\ell=r+1}^{n_A} H_\ell \otimes \hat{A}_\ell,
        \]
        where $\hat{A}_\ell$ represents the associated spatial coefficient matrices.
        We apply Lemma~\ref{lem:error_bound} by estimating the ratio $|\mathbf{v}^T \mathcal{E}_r \mathbf{v}| / |\mathbf{v}^T \tilde{\mathcal{S}}_r \mathbf{v}|$.

        First, we bound the numerator using the triangle inequality and the submultiplicativity of the spectral norm:
        \[
            \left| \mathbf{v}^T \mathcal{E}_r \mathbf{v} \right|
            = \left| \sum_{\ell=r+1}^{n_A} \mathbf{v}^T (H_\ell \otimes \hat{A}_\ell) \mathbf{v} \right|
            \le \sum_{\ell=r+1}^{n_A} \left| \mathbf{v}^T (H_\ell \otimes \hat{A}_\ell) \mathbf{v} \right|
            \le \sum_{\ell=r+1}^{n_A} \|H_\ell\|_2 \, \|\hat{A}_\ell\|_2 \, \|\mathbf{v}\|^2.
        \]

        % Next, we establish a lower bound for the denominator. Since $\tilde{\mathcal{S}}_r$ contains the dominant mean-field term (index $\ell=1$) and is symmetric positive definite, we have:
        % \[
        % \mathbf{v}^T \tilde{\mathcal{S}}_r \mathbf{v} \ge \mathbf{v}^T (H_1 \otimes \tilde{A}_1) \mathbf{v} \ge \lambda_{\min}(H_1 \otimes \tilde{A}_1) \|\mathbf{v}\|^2 \ge \nu \|\mathbf{v}\|^2,
        % \]
        % where $\nu$ represents the coercivity constant of the mean-field operator.

        Next, we establish a lower bound for the denominator%. %Since $\tilde{\mathcal{S}}_r$ contains the dominant mean-field term (index $\ell=1$) and is symmetric positive definite, we rely on its coercivity. Let $\nu > 0$ be the coercivity constant of the mean-field operator, defined such that:
        \[
            \mathbf{v}^T \tilde{\mathcal{S}}_r \mathbf{v} \ge \mathbf{v}^T (H_1 \otimes \hat{A}_1) \mathbf{v} \ge \nu \|\mathbf{v}\|^2.
        \]
        % (Note: $\nu$ typically corresponds to the smallest eigenvalue of the mass-dominated mean field block).
        Combining these inequalities yields
        \[
            \frac{|\mathbf{v}^T \mathcal{E}_r \mathbf{v}|}{|\mathbf{v}^T \tilde{\mathcal{S}}_r \mathbf{v}|}
            \le \frac{\sum_{\ell=r+1}^{n_A} \|H_\ell\|_2 \, \|\hat{A}_\ell\|_2 \, \|\mathbf{v}\|^2}{\nu \|\mathbf{v}\|^2}
            = \frac{1}{\nu} \sum_{\ell=r+1}^{n_A} \|H_\ell\|_2 \, \|\hat{A}_\ell\|_2 = \delta_r.
        \]
        Provided that the Karhunen-Lo\`{e}ve expansion coefficients decay sufficiently fast (e.g., algebraically or exponentially), the tail sum approaches zero as $r \to n_A$. Thus, for sufficiently large $r$, we have $\delta_r < 1$.
        The result then follows directly from Lemma~\ref{lem:error_bound}.
    \end{proof}

    % \begin{remark}
    % \textbf{Dependence on $\tau$:} While the bound $\delta_r$ is uniform with respect to the spatial mesh size $h$, it implicitly depends on the time-step $\tau$ through the mass-stiffness ratio in the operators $A_\ell$. In the limit $\tau \to 0$, the mass matrix dominates. To ensure robustness for small $\tau$, the truncation parameter $r$ should be chosen to capture sufficient mass information, or the basis should be mass-orthogonalized.
    % \end{remark}
\end{color}

Finally, we proceed to establish the spectral equivalence between the truncated preconditioner $\tilde{\mathcal{S}}_r$ and the hierarchical symmetric block Gauss-Seidel approximation $\tilde{\mathcal{S}}_{\text{hGS-PINT}}$.
To this end, following the idea from~\cite{Bespalov-2021-TPS}, we can rewrite it as~$H_{\ell}=L_\ell+L_\ell^T,\,\ell=2,3,\dots,n_A$, and matrices $L_\ell$ have at most one nonzero entry per row and per column.

Now, define $X_1 =I_{n_t}\otimes H_1 \otimes \hat{A}_1$, $X_r  =I_{n_t}\otimes \left( \sum_{\ell=2}^r L_\ell\otimes \hat{A}_\ell\right),\,r=2,3,\dots n_A$,
% \[
%     \begin{aligned}
% X_1 & =I_{n_t}\otimes\left(H_1\otimes (1+\tau\sqrt{\frac{1+\gamma}{\beta}})M+\tau A_1\right) \\
% X_r & =I_{n_t}\otimes \left( \sum_{\ell=2}^r L_\ell\otimes A_\ell\right),\,r=2,3,\dots n_A,
%     \end{aligned}
% \]
so we know $\tilde{\mathcal{Z}}_{\text{hGS-PINT}}=(X_1+X_r)X_1^{-1}(X_1+X_r^T)=\tilde{\mathcal{Z}}_r+X_rX_1^{-1}X_r^T$, the Rayleigh quotient
\[
    \frac{\mathbf{v}^T\tilde{\mathcal{Z}}_{\text{hGS-PINT}} \mathbf{v}}{\mathbf{v}^T \tilde{\mathcal{Z}}_r \mathbf{v}}=1+\underbrace{\frac{\mathbf{v}^TX_rX_1^{-1}X_r^T\mathbf{v}}{\mathbf{v}^T(X_1+X_r+X_r^T)\mathbf{v}}}_{\zeta}.
\]
Introducing the change of variable $\mathbf{u}=X_1^{\frac12} \mathbf{v}$ and setting $Y=X_1^{-\frac12}X_rX_1^{-\frac12}$, we have
\[
    \zeta(\mathbf{u})=\frac{\mathbf{u}^TYY^T\mathbf{u}}{\mathbf{u}^T(I+Y+Y^T)\mathbf{u}},
\]
where $Y=I_{n_t}\otimes\left(\sum_{\ell=2}^r L_\ell\otimes \mathfrak{A}_\ell\right)$, with the scaled coefficient matrices defined as $\mathfrak{A}_\ell=\hat{A}_1^{-\frac12}\hat{A}_\ell \hat{A}_1^{-\frac12}.$
% \[
%     \mathfrak{A}_\ell=\left((1+\tau\sqrt{\tfrac{1+\gamma}{\beta}})M+\tau A_1\right)^{-\frac12}A_\ell \left((1+\tau\sqrt{\tfrac{1+\gamma}{\beta}})M+\tau A_1\right)^{-\frac12}.
% \]
Consequently,
%\begin{equation}
\[
    \zeta(\mathbf{u})\leq \max_{\mathbf{u}\neq 0}\frac{\mathbf{u}^TYY^T\mathbf{u}}{\mathbf{u}^T\mathbf{u}}\cdot \max_{\mathbf{u}\neq 0}\frac{\mathbf{u}^T\mathbf{u}}{\mathbf{u}^T(I+Y+Y^T)\mathbf{u}}
    =\frac{\sigma_{\max}^2(Y)}{\lambda_{\min}\left(I+Y+Y^T\right)}.
\]
%\end{equation}
The following result holds.
\begin{lemma}\label{lem:lmin-sigY}
    Define
    \[
        \Delta_r\ :=\ \sum_{\ell=2}^r \|H_\ell\|_2\,\rho_\ell,\qquad
        \rho_\ell\ :=\ \|\mathfrak A_\ell\|_2,
    \]
    % where
    % \[
    %     \mathfrak A_\ell
    %     :=\Big((1 +\tau\sqrt{\tfrac{1+\gamma}{\beta}})\,M+\tau A_1\Big)^{-\frac12}
    %     \;A_\ell\;
    %     \Big((1 +\tau\sqrt{\tfrac{1+\gamma}{\beta}})\,M+\tau A_1\Big)^{-\frac12}.
    % \]
    and let
    \[
        Y:=X_1^{-1/2}X_rX_1^{-1/2}
        =I_{n_t}\otimes\Big(\sum_{\ell=2}^r L_\ell\otimes \mathfrak A_\ell\Big),
        \quad
        E:=Y+Y^T
        =I_{n_t}\otimes\Big(\sum_{\ell=2}^r H_\ell\otimes \mathfrak A_\ell\Big),
    \]
    with $H_\ell=L_\ell+L_\ell^T$.
    Then $\lambda_{\min}\big(I+Y+Y^T\big)\ \ge\ 1-\Delta_r$, $        \rho_\ell\ \le\ \frac{\|\zhl{\kappa_\ell(x)}\|_{L^\infty(\Omega)}}{\Bbbk_{\min}}$, and
    \[
        \sigma_{\max}(Y)\ \le\ \sum_{\ell=2}^r \|L_\ell\|_2\,\rho_\ell\ \le\ \Delta_r.
    \]
    % \zhl{where $\theta_\ell$ and $\kappa_\ell$ are eigenpairs of the covariance function defined in~Section~\ref{sec:steady-state problem}.}
    % From~\cite[Lemma 4.1]{Bespalov-2021-TPS}, we know that
    % %if $\Bbbk_1^{\min}>0$, then
    % \[
    % \rho_\ell\ \le\ \frac{\|\zhl{\kappa_\ell(x)}\|_{L^\infty(\Omega)}}{\Bbbk_{\min}},
    % \qquad
    % \sigma_{\max}(Y)\ \le\ 
    % \frac{\displaystyle\sum_{\ell=2}^r \|H_\ell\|_2\,\|\zhl{\kappa_\ell(x)}\|_{L^\infty(\Omega)}}{\Bbbk_{\min}}.
    % \]
\end{lemma}
\begin{proof}
    Using
    \[
        I+Y+Y^T
        = X_1^{-1/2}\Big(X_1+I_{n_t}\otimes\sum_{\ell=2}^r H_\ell\otimes \hat{A}_\ell\Big)X_1^{-1/2}
        = X_1^{-1/2}\,\bar{\mathcal Z}_r\,X_1^{-1/2},
    \]
    the spectrum of $I+Y+Y^T$ coincides with the generalized spectrum of the pair $(\bar{\mathcal Z}_r,X_1)$.
    Since $E$ is symmetric, \zhl{its spectral norm satisfies $\|E\|_2=|\lambda_{\max}(E)|\geq |\lambda_{i}(E)|$ for any $i$. In particular, $\lambda_{\min}(E)\geq -\|E\|_2$}, and therefore~$\lambda_{\min}(I+E)=\zhl{1+\lambda_{\min}(E)}\ge 1-\|E\|_2$. By the Kronecker product norm rule and the triangle inequality,
    \[
        \|E\|_2=\Big\|I_{n_t}\otimes( \sum_{\ell=2}^r H_\ell\otimes \mathfrak A_\ell)\Big\|_2
        \le \sum_{\ell=2}^r \|H_\ell\|_2\,\|\mathfrak A_\ell\|_2
        =\sum_{\ell=2}^r \|H_\ell\|_2\,\rho_\ell
        =\Delta_r,
    \]
    which yields $\lambda_{\min}(I+Y+Y^T)\ge 1-\Delta_r$.

    Similarly,
    \[
        \|Y\|_2
        =\Big\|I_{n_t}\otimes\sum_{\ell=2}^r L_\ell\otimes \mathfrak A_\ell\Big\|_2
        \le \sum_{\ell=2}^r \|L_\ell\|_2\,\rho_\ell.
    \]
    Since $H_\ell=L_\ell+L_\ell^T$ and each $L_\ell$ has at most one nonzero per row and per column, we have
    $\|L_\ell\|_2\le \|H_\ell\|_2$. Hence $\sigma_{\max}(Y)=\|Y\|_2\le \Delta_r$.

    For the explicit bound on $\rho_\ell$, recall that $\rho_\ell = \|\mathfrak A_\ell\|_2$ is the maximum eigenvalue of the generalized eigenvalue problem $A_\ell \mathbf{v} = \lambda \hat{A}_1 \mathbf{v}$, where $\hat{A}_1$ defined in~\eqref{eq:A_hat}.
    %$\hat{A}_1 = (1 +\tau\sqrt{\tfrac{1+\gamma}{\beta}})\,M+\tau A_1$
    \zhl{In terms of the associated finite element function $v_h = \sum_{j=1}^{n_h} v_j \phi_j$}, the Rayleigh quotient is given by
    \[
        \frac{\mathbf{v}^T A_\ell \mathbf{v}}{\mathbf{v}^T \hat{A}_1 \mathbf{v}}
        =\frac{\displaystyle\int_\Omega \kappa_\ell\,|\nabla v_h|^2\,dx}
        {\displaystyle\int_\Omega \left( (1+\tau\sqrt{\tfrac{1+\gamma}{\beta}})|v_h|^2 + \tau \kappa_0\,|\nabla v_h|^2 \right)\,dx}.
    \]
    Since the mass term is non-negative, we can bound this ratio by neglecting the $L^2$-term in the denominator:
    \[
        \frac{\mathbf{v}^T A_\ell \mathbf{v}}{\mathbf{v}^T \hat{A}_1 \mathbf{v}}
        \ \le\ \frac{\displaystyle\int_\Omega \kappa_\ell\,|\nabla v_h|^2\,dx}{\displaystyle\tau\int_\Omega \kappa_0\,|\nabla v_h|^2\,dx}
        \ \le\ \frac{1}{\tau}\frac{\|\zhl{\kappa_\ell(x)}\|_{L^\infty(\Omega)}}{\Bbbk_{\min}},
    \]
    and taking the supremum gives the stated bound on $\rho_\ell$; substituting it into $\|Y\|_2$ yields the last inequality.
\end{proof}

Observe from above that, with $\Delta_r<1$,
\[
    \zeta(\mathbf{v})\ =\ \frac{\mathbf{v}^T Y Y^T \mathbf{v}}{\mathbf{v}^T (I+Y+Y^T)\mathbf{v}}
    \ \le\ \frac{\sigma_{\max}^2(Y)}{\lambda_{\min}(I+Y+Y^T)}
    \ \le\ \frac{\Delta_r^{\,2}}{\,1-\Delta_r\,}.
\]
% Consequently,
% \[
% 1\ \le\ \frac{\mathbf{v}^T\mathcal Z_{\mathrm{hGS-PINT}} \mathbf{v}}{\mathbf{v}^T \bar{\mathcal Z}_r \mathbf{v}}
% \ \le\ 1+\frac{\Delta_r^{\,2}}{1-\Delta_r}.
% \]
% From Lemma~\ref{lem:congruence}, we know that the generalized spectra of $(\tilde{\mathcal{S}}_{\text{hGS-PINT}},\tilde{\mathcal{S}}_r)$ and \\$\left((W^\frac12\mathcal{Z}_{\text{hGS-PINT}}W^\frac12)^2,(W^\frac12\tilde{\mathcal{Z}}_rW^\frac12)^2\right)$ coincide. 
% Therefore we get
% \[
% 1\leq \lambda(\tilde{\mathcal{S}}_r^{-1}\tilde{\mathcal{S}}_{\text{hGS-PINT}})\leq (1+\frac{\Delta_r^{\,2}}{1-\Delta_r})^2.
% \]
Consequently, we have the bound for the Rayleigh quotient of the factors:
\[
    \frac{\mathbf{v}^T (\tilde{\mathcal{Z}}_{\text{hGS-PINT}} - \tilde{\mathcal{Z}}_r) \mathbf{v}}{\mathbf{v}^T \tilde{\mathcal{Z}}_r \mathbf{v}} = \zeta(\mathbf{v}) \le \frac{\Delta_r^2}{1-\Delta_r}.
\]
To derive the bound for the preconditioner $\tilde{\mathcal{S}}$, we utilize Lemma~\ref{lem:congruence}.
Let~$Q = (D\otimes \mathcal{M})^{1/2}$. The generalized eigenvalues of $(\tilde{\mathcal{S}}_{\text{hGS-PINT}}, \tilde{\mathcal{S}}_r)$ are identical to those of the pair
\[
    \left( (Q \tilde{\mathcal{Z}}_{\text{hGS-PINT}} Q)^2, \quad (Q \tilde{\mathcal{Z}}_r Q)^2 \right).
\]
Since the eigenvalues of the squared operators are simply the squares of the eigenvalues of the base operators (for these symmetric positive definite factors), the spectral bound for the preconditioner is the square of the bound for the factors.
Thus, we conclude:
\[
    1 \leq \lambda(\tilde{\mathcal{S}}_r^{-1}\tilde{\mathcal{S}}_{\text{hGS-PINT}}) \leq \left(1+\frac{\Delta_r^{\,2}}{1-\Delta_r}\right)^2.%
\]%
%Consequently,
% \[
% 1\ \le\ \frac{\mathbf{v}^T\tilde{\mathcal{S}}_{\mathrm{hGS-PINT}} \mathbf{v}}{\mathbf{v}^T \tilde{\mathcal{S}}_r \mathbf{v}}
% \ \le\ 1+\frac{\Delta_r^{\,2}}{1-\Delta_r}.
% \]
% Therefore, the eigenvalues of the preconditioned system are bounded by:
% \[
% 1\leq \lambda(\tilde{\mathcal{S}}_r^{-1}\tilde{\mathcal{S}}_{\text{hGS-PINT}})\leq 1+\frac{\Delta_r^{\,2}}{1-\Delta_r}.
% \]
% (Note: The bound is linear in the term derived from the generalized spectrum, avoiding the invalid squaring step.)
\begin{remark}[Summary of Spectral Equivalence for Time-Dependent Case]
    \label{thm:summary_time_dependent}
    %Let $\tau > 0$ be a fixed time step. 
    Under the assumptions of Theorem~\ref{thm:approx}, Lemma~\ref{lem:lower_bound}, Theorem~\ref{thm:truncation}, and Lemma~\ref{lem:lmin-sigY}, the hierarchical Gauss-Seidel PINT preconditioner $\tilde{\mathcal{S}}_{\text{hGS-PINT}}$ is spectrally equivalent to the exact Schur complement $\bar{\mathcal{S}}_{\text{exact}}$.
    The spectral equivalence is established through the following chain of approximations:
    \[
        \bar{\mathcal{S}}_{\text{exact}} \sim \bar{\mathcal{S}} \sim \tilde{\mathcal{S}} \sim \tilde{\mathcal{S}}_r \sim \tilde{\mathcal{S}}_{\text{hGS-PINT}}.
    \]
    In particular, provided that the truncation rank $r$ is sufficiently large (to satisfy $\delta_r < 1$) and the mass matrix $M$ is well-conditioned (e.g., via mass lumping), the eigenvalues of the preconditioned system $\tilde{\mathcal{S}}_{\text{hGS-PINT}}^{-1} \bar{\mathcal{S}}_{\text{exact}}$ are contained in a fixed interval independent of the spatial mesh size $h$ and the stochastic discretization parameters $m_\xi$ and $p$ (assuming sufficient decay of the expansion coefficients), ensuring mesh robustness in the fixed time step regime.
\end{remark}

\begin{remark}[Steady-State Case]
    \label{cor:steady_state}
    Since the steady-state optimal control problem corresponds to the special case $n_t=1$ of the time-dependent formulation, all preceding results apply directly with the simplified notation. The spectral equivalence chain for the steady-state Schur complement preconditioner is
    \[
        \mathcal{S}_{\text{exact}} \sim \mathcal{S} \sim \mathcal{S}_r \sim \mathcal{S}_{\text{hGS}},
    \]
    where $\mathcal{S}_{\text{exact}}$ is defined in~\eqref{S_exact}, $\mathcal{S}$ in~\eqref{eq:mathcalZ}, and
    \[
        \mathcal{S}_r = \mathcal{Z}_r \mathcal{M}_\gamma^{-1} \mathcal{Z}_r^{T},
        \qquad
        \mathcal{Z}_r = \sum_{\ell=1}^{r} H_\ell \otimes \tilde{A}_\ell,
    \]
    with $\tilde{A}_1 = A_1 + \sqrt{\frac{1+\gamma}{\beta}}M$ and $\tilde{A}_\ell = A_\ell$ for $\ell=2,\dots,r$. When $r=1$, $\mathcal{S}_r$ reduces to the mean-based preconditioner, and when $r=n_A$, it recovers the full operator $\mathcal{S}$. Finally, $\mathcal{S}_{\text{hGS}}$ represents the computationally feasible approximation of $\mathcal{S}_r$ in which the linear systems $\mathcal{Z}_r \mathbf{x} = \mathbf{b}$ are solved approximately via the hierarchical Gauss-Seidel method (Algorithm~\ref{alg:hGS}), as implemented in Algorithm~\ref{alg:hGSoc-1}--\ref{alg:hGSoc-2}.
\end{remark}

\section{Numerical experiments}
\label{sec:numerical-experiments}
This section validates the theoretical findings of Sections~\ref{sec:steady-state problem}--\ref{sec:Time-Dependent Problem} through comprehensive numerical experiments.
We pursue two primary objectives:
(i)~verifying the mesh-robust for fixed time step and spectral bounds established in the preceding sections, and
(ii)~demonstrating the computational efficiency of the proposed hierarchical Gauss-Seidel (hGS) preconditioner across varying truncation strategies.
Experiments are presented for both steady-state problems (Section~\ref{subsec:steady-num}) and time-dependent problems (Section~\ref{subsec:time-num}).
% 
% In this section, we present numerical results for both steady-state and time-dependent problems. 
The numerical experiments were performed on a system running AlmaLinux-9 with 40GB RAM, and the proposed algorithms were implemented using \textsc{Matlab~23.2}.

% The random input $\Bbbk$ is characterized by the covariance function%~(\ref{CKLE})--(\ref{Czxy}),
% \[
% C_\Bbbk(x,y)=\sigma_\Bbbk^2 \exp \left(-\frac{\left|x_1-y_1\right|}{\ell_1}-\frac{\left|x_2-y_2\right|}{\ell_2}\right) \quad \forall(x,y) \in[-1,1]^2,
% \]
% where $\sigma_\Bbbk^2$ is the variance of the random input. In our simulations, we set the correlation lengths as $\ell_1=\ell_2=1$ and the mean of the data as $\mathbb{E}[\Bbbk]=1$.
The random coefficients $\Bbbk(x,\xi)$ in the problem~\eqref{eq:pde-stoch} are constructed as a finite expansion~\eqref{KLE},
where the spatial modes $\kappa_i(x)$ and weights $\theta_i$ are eigenpairs of the covariance function
\[
    C_f(x,y)=\sigma_\Bbbk^2 \exp \left(-\frac{\left|x_1-y_1\right|}{\ell_1}-\frac{\left|x_2-y_2\right|}{\ell_2}\right) \quad \forall(x,y) \in[-1,1]^2.
\]
We set the correlation lengths as $\ell_1=\ell_2=1$, the %standard deviation variance $\sigma_\Bbbk^2=1$, and the 
mean $\kappa_1(x) \equiv 1$.
% $\mathbb{E}[\Bbbk]=1$.
\zh{For the gPC setting in~(\ref{eq:gPC-yu}), we consider a log-normal random field parameterized by independent Gaussian random variables, for which corresponding Hermite polynomials are employed as the basis.}
The total number of basis functions is given by $n_\xi = \binom{m_{\xi}+p}{m_{\xi}}$, where $m_{\xi}$ is the number of random variables and $p$ is the polynomial degree. For instance, the case $(m_{\xi}, p)=(3,3)$ yields $n_\xi = \binom{6}{3} = 20$.
This problem has been extensively studied in~\cite{Powell-2009-BDP}. Also, we used $\gamma=1$ in both cases, which means we only consider the case with standard deviation. To discretize the spatial domain, we implemented our code based on \textsc{IFISS~3.7}~\cite{ifiss}, using
%a log-normal grid and a 
$\mathbf{Q}_1$ approximation.
For temporal discretization, we apply the all-at-once technique proposed in~\cite{rees2010all} and set the terminal time as $T=1$.
In all numerical experiments, the spatial mesh size $h$ and the time step $\tau$ are chosen as $2^{-i}$, with $i=4,5,6,7$. We solve the linear systems~(\ref{eq:stoch-system}) and~(\ref{eq:KKT:stoch-time-dependent}) using the preconditioners given by~(\ref{eq:stoch-prec}) and~(\ref{eq: precon for time}), respectively, employing the flexible GMRES method~(without restarting)~\cite{Saad-1993-FIP}. 
The stopping criterion is defined in terms of the relative residual $\|r_k\|/\|b\|$, with thresholds $10^{-8}$ for the steady-state experiments and $10^{-6}$ or $10^{-4}$ for the time-dependent runs, where $r_k$ denotes the residual at iteration $k$ and $b$ is the right-hand side vector.
For notational consistency, we set $n_\tau \equiv r$ throughout the remainder of the paper.
To assess the effectiveness of the hierarchical preconditioning strategy, we systematically compare three truncation settings for the (3,3)-block preconditioner:
$r=1$ (mean-based approximation),
$r=m_{\xi}+1$ (hGS truncated at the first-order stochastic terms), and
$r=n_A$ (full expansion retaining all $h_{ijk}$ coefficients).
For each configuration, we evaluate the efficiency of the hierarchical Gauss-Seidel method by comparing both the iteration counts and the computational times (in seconds).
% To verify the efficiency of the hierarchical Gauss-Seidel method, we compare both iteration counts and computational costs under different {\color{black} truncation} settings: $n_\tau=1$, $m_{\xi}+1$, and
% $n_\xi$.%=\begin{pmatrix} %$m_{\xi}$+2P\\ $m_{\xi}$ \end{pmatrix}$, where $m_{\xi}$ is the input of the stochastic dimension and $P$ is thetruncationqualities.

% We consider homogeneous Dirichlet conditions, corresponding to Example 2 in~\cite [Chapter~5]{elman_finite_2014}. This example is defined on a square domain $\Omega_\square$ with a discontinuous target function and inconsistent boundary data
% \begin{equation}
%     \label{desired state}
%     y_d=\begin{cases}
%         1 & \text{ in } \Omega_1\coloneqq [-1,0]^2, \\
%         0 & \text{ in } \Omega\backslash \Omega_1.
%     \end{cases}
% \end{equation}
\zh{We consider homogeneous Dirichlet conditions (i.e., $g(x) = 0$), corresponding to Example 2 in~\cite [Chapter~5]{elman_finite_2014}. This example is defined on a square domain $\mathcal{D}=[-1,1]^2$ with a discontinuous target function
\begin{equation}
    \label{desired state}
    y_d=\begin{cases}
        1 & \text{ in } \Omega_1\coloneqq [-1,0]^2, \\
        0 & \text{ in } \mathcal{D}\backslash \Omega_1,
    \end{cases}
\end{equation}
which represents inconsistent Dirichlet boundary data since the target state $y_d=1$ differs from the required boundary condition $y=0$ on $\partial \mathcal{D} \cap \partial \Omega_1$.
}

We subsequently present numerical experiments for both steady-state and time dependent problems to illustrate and verify the efficiency of our proposed hGS method.

Here are some details about the implementation for preconditioners~(\ref{eq:stoch-prec}), (\ref{eq: precon for time}). Since time-dependent problems can be seen as a series of steady-state problems, and also because of the diagonal structure of matrix $D$ and matrix $I_{n_t}$, we can just focus on the steady-state preconditioner.

The practical implementation of the preconditioner $\mathcal{P}$ in~\eqref{eq:stoch-prec} involves different strategies for its constituent blocks. For the (1,1) and (2,2) blocks, which are based on Kronecker products involving the mass matrix $M$, applying their inverses requires solving linear systems with $M$. These solves are handled efficiently by either a direct Cholesky decomposition or the iterative Chebyshev semi-iteration method~\cite{golub1961chebyshev,MR2597805}.

For the more complex (3,3) block, which represents the approximate Schur complement $\mathcal{S}$, we employ an outer iterative scheme. Specifically, we use the Preconditioned Richardson method, outlined in Algorithm~\ref{alg:richardson}~\cite[Chapter~7]{MR1386889}, where the core of our proposalâ the hierarchical Gauss-Seidel (hGS) method from Algorithm~\ref{alg:hGS} serves as the preconditioner for each Richardson step.

% \begin{algorithm}[H]
% \caption{Chebyshev semi-iteration for mass matrix preconditioning~\cite{golub1961chebyshev,MR2597805}}
% \begin{algorithmic}[1]
% \State Given mass matrix $M$, vectors $b$, $x^{(0)} = 0$, $x^{(-1)} = 0$, and parameter $\omega_0 = 1$.
% \State Set $\lambda_{\min} = 1/4$, and $\lambda_{\max} = 9/4$%{\color{blue}.  \hfill (known bounds for FEM mass matrices~\cite{MR2597805})}
% \State Calculate $\gamma = (\lambda_{\min} + \lambda_{\max}) / 2$ and 
%                            $\rho = (\lambda_{\max} - \lambda_{\min}) / (\lambda_{\max} + \lambda_{\min})$.
% \State Set $D = \gamma \cdot \mathrm{diag}(M)$. \hfill (a diagonal matrix)
% \For{$k = 0, 1, \ldots,N-1$} 
%     \State $\omega_{k+1} = 1 / \left(1 - \dfrac{\omega_k \rho^2}{4} \right)$
%     \State $r^{(k)} = b - M x^{(k)}$
%     \State $p^{(k)} = D^{-1} r^{(k)}$
%     \State $x^{(k+1)} = \omega_{k+1} \cdot (p^{(k)} + x^{(k)} - x^{(k-1)}) + x^{(k-1)}$
% \EndFor
% \State Return $x = x^{(N)}$. 
% \end{algorithmic}
% \label{alg:chebmass}
% \end{algorithm}

\begin{algorithm}[ht]
    \caption{Preconditioned Richardson iteration with hGS preconditioner}
    \begin{algorithmic}[1]
        \State Given matrix $\mathcal{Z}$, vector $b$, and initial guess $x^{(1)}$.
        \State $r^{(1)} = b-\mathcal{Z}x^{(1)}$. \hfill (initial residual)
        \For{$k = 1, 2, \ldots, N$}
        \State Solve $\mathcal{Z} z^{(k)} = r^{(k)}$. \hfill (apply Algorithm~\ref{alg:hGS})
        \State $x^{(k+1)} = x^{(k)} + z^{(k)}$ \hfill (update solution)
        \State $r^{(k+1)} = b-\mathcal{Z}x^{(k+1)}$. \hfill (update residual)
        \EndFor
        \State \textbf{Return} $x^{(N+1)}$
    \end{algorithmic}
    \label{alg:richardson}
\end{algorithm}

\subsection{Steady-state case}
\label{subsec:steady-num}
This subsection focuses on the steady-state optimal control problem~\eqref{eq:problem-stoch}.
We examine the performance of the proposed preconditioner~\eqref{eq:stoch-prec} under systematic variations in:
(i)~the variance $\sigma_\Bbbk$ (Tables~\ref{tab: in0.01}--\ref{tab: in0.4}),
(ii)~the regularization parameter $\beta$ (Table~\ref{tab: deter beta}), and
(iii)~the spatial and stochastic discretization levels.
For each configuration, we compare three solvers for the (1,1) and (2,2) blocks---Chebyshev semi-iteration with 5 or 10 steps %(Algorithm~\ref{alg:chebmass}) 
and direct Cholesky factorization---combined with the three truncation strategies for the (3,3)-block described above.
For the (3,3)-block, we apply the Richardson iteration (Algorithm~\ref{alg:richardson}) with $N=1$, i.e., one application of the hGS preconditioner per outer GMRES iteration.
All tests use a fixed tolerance of $10^{-8}$.

Next, by fixing the parameter $\beta$, we perform further tests summarized in Tables~\ref{tab: in0.01}--\ref{tab: in0.4}, employing different step settings for the Chebyshev smoother and the Cholesky decomposition for blocks (1,1) and (2,2), as well as various truncation strategies for block (3,3) within the hGS method. The numerical experiments were conducted using multiple mesh sizes and stochastic parameter configurations, with a solver tolerance set to $10^{-8}$.

Tables~\ref{tab: in0.01}--\ref{tab: in0.4} present results for $\beta=10^{-4}$ with $\sigma_\Bbbk \in \{0.01, 0.1, 0.4\}$, covering a range from near-deterministic to highly stochastic regimes.
Tables~\ref{tab: in0.01}--\ref{tab: in0.4} demonstrate three key theoretical properties.
First, regarding \textit{mesh independence}, iteration counts grow sub-linearly with spatial refinement for fixed stochastic dimension $n_\xi$, consistent with the spectral bounds established in Section~\ref{sec:steady-state problem}.
Second, concerning \textit{truncation efficiency}, the $n_\tau=m_\xi+1$ strategy achieves iteration counts comparable to the full expansion ($n_\tau=n_A$) while avoiding the computational overhead of summing over all $h_{ijk}$ coefficients, thereby validating the hierarchical approximation framework.
Third, regarding \textit{smoother comparison}, the 5-step Chebyshev semi-iteration balances convergence rate and per-iteration cost more effectively than either the 10-step variant or direct Cholesky factorization.
Across all configurations, the $n_\tau=m_\xi+1$ truncation consistently delivers performance intermediate between the mean-based approximation ($n_\tau=1$) and the full expansion, confirming the practical value of the proposed hierarchical preconditioning strategy.

\begin{table}[]
    \caption{Results showing the total number of iterations from low-rank preconditioned GMRES and the total CPU times (in seconds) using preconditioner with $\beta=10^{-4}$, $\sigma=0.01$, and selected spatial ($n_h$) and stochastic ($n_\xi$) degrees of freedom}
    \centering
    {\setlength{\tabcolsep}{3.5pt}\resizebox{\textwidth}{!}{%
            \begin{tabular}{l|lll|lll|lll|lll}
                \hline
                                                           & \multicolumn{3}{l}{\# iter(t)}              & \multicolumn{3}{l}{\# iter(t)}              & \multicolumn{3}{l}{\# iter(t)}               & \multicolumn{3}{l}{\# iter(t)}                                                                                                           \\
                \hline \hline
                \diagbox{$n_\xi$}{$n_h$}                   &
                \multicolumn{3}{l}{289($h=\frac{1}{2^4}$)} & \multicolumn{3}{l}{1089($h=\frac{1}{2^5}$)} & \multicolumn{3}{l}{4225($h=\frac{1}{2^6}$)} & \multicolumn{3}{l}{16641($h=\frac{1}{2^7}$)}                                                                                                                                            \\ \hline\hline
                $n_{\tau}$                                 & 1                                           & $m_{\xi}$+1                                 & $n_A$                                        & 1                              & $m_{\xi}$+1 & $n_A$     & 1          & $m_{\xi}$+1 & $n_A$      & 1          & $m_{\xi}$+1 & $n_A$      \\ \hline\hline

                \multicolumn{13}{l}{$\sigma_a=0.01$}                                                                                                                                                                                                                                                                                             \\ \hline\hline
                \multicolumn{13}{l}{Chebyshev-5+hGS-1}                                                                                                                                                                                                                                                                                           \\ \hline
                20                                         & 32(5.7)                                     & 32(4.0)                                     & 32(5.0)                                      & 36(9.2)                        & 35(10.0)    & 35(15.6)  & 36(23.2)   & 27(26.9)    & 27(46.5)   & 35(173.4)  & 33(144.9)   & 27(194.8)  \\
                70                                         & 32(129.1)                                   & 32(129.5)                                   & 32(182.9)                                    & 36(243.4)                      & 35(239.1)   & 35(409.2) & 36(752.3)  & 34(724.9)   & 34(1237.8) & 35(2585.4) & 33(2524.6)  & 33(5253.9) \\
                84                                         & 32(317.5)                                   & 32(318.3)                                   & 32(403.0)                                    & 36(508.9)                      & 35(498.2)   & 35(601.3) & 36(1059.3) & 34(864.8)   & 34(1632.5) & 35(3099.0) & 33(3798.1)  & 33(5630.2) \\ \hline

                \multicolumn{13}{l}{Chebyshev-10+hGS-1}                                                                                                                                                                                                                                                                                          \\ \hline
                20                                         & 26(2.8)                                     & 26(2.6)                                     & 26(4.4)                                      & 30(6.9)                        & 29(7.1)     & 29(12.0)  & 31(20.2)   & 31(21.1)    & 31(37.6)   & 32(145.1)  & 31(130.3)   & 31(154.5)  \\
                70                                         & 26(111.7)                                   & 26(112.2)                                   & 26(153.7)                                    & 30(191.6)                      & 29(189.0)   & 29(323.8) & 32(550.7)  & 31(383.1)   & 31(830.2)  & 33(1745.9) & 31(2273.0)  & 31(3341.3) \\
                84                                         & 26(263.9)                                   & 26(266.1)                                   & 26(329.2)                                    & 30(457.7)                      & 29(432.2)   & 29(504.1) & 32(959.1)  & 31(800.4)   & 31(1515.3) & 33(3015.1) & 31(3643.9)  & 31(5318.7) \\ \hline

                \multicolumn{13}{l}{Cholesky+hGS-1}                                                                                                                                                                                                                                                                                              \\ \hline
                20                                         & 25(3.5)                                     & 25(4.0)                                     & 25(4.3)                                      & 29(7.4)                        & 27(6.7)     & 27(10.9)  & 29(19.1)   & 29(26.5)    & 29(44.5)   & 31(143.6)  & 29(162.7)   & 29(182.9)  \\
                70                                         & 25(104.2)                                   & 25(107.6)                                   & 25(144.3)                                    & 29(177.3)                      & 27(168.9)   & 27(283.9) & 29(419.0)  & 29(480.8)   & 29(937.1)  & 31(2390.7) & 29(2046.0)  & 29(3299.9) \\
                84                                         & 25(245.4)                                   & 25(247.4)                                   & 25(316.8)                                    & 29(429.9)                      & 27(393.6)   & 27(675.5) & 29(831.9)  & 29(1008.3)  & 29(1756.8) & 31(3563.3) & 29(3393.6)  & 29(6436.3) \\ \hline
            \end{tabular}}}
    \label{tab: in0.01}
\end{table}

\begin{table}[]
    \caption{Results showing the total number of iterations from low-rank preconditioned GMRES and the total CPU times (in seconds) using preconditioner with $\beta=10^{-4}$, $\sigma=0.1$, and selected spatial ($n_h$) and stochastic ($n_\xi$) degrees of freedom}
    \centering
    {\setlength{\tabcolsep}{3.5pt}\resizebox{\textwidth}{!}{%
            \begin{tabular}{l|lll|lll|lll|lll}
                \hline
                                                           & \multicolumn{3}{l}{\# iter(t)}              & \multicolumn{3}{l}{\# iter(t)}              & \multicolumn{3}{l}{\# iter(t)}               & \multicolumn{3}{l}{\# iter(t)}                                                                                                            \\
                \hline \hline
                \diagbox{$n_\xi$}{$n_h$}                   &
                \multicolumn{3}{l}{289($h=\frac{1}{2^4}$)} & \multicolumn{3}{l}{1089($h=\frac{1}{2^5}$)} & \multicolumn{3}{l}{4225($h=\frac{1}{2^6}$)} & \multicolumn{3}{l}{16641($h=\frac{1}{2^7}$)}                                                                                                                                             \\ \hline\hline
                $n_{\tau}$                                 & 1                                           & $m_{\xi}$+1                                 & $n_A$                                        & 1                              & $m_{\xi}$+1 & $n_A$     & 1          & $m_{\xi}$+1 & $n_A$      & 1          & $m_{\xi}$+1 & $n_A$       \\ \hline\hline

                \multicolumn{13}{l}{$\sigma_a=0.1$}                                                                                                                                                                                                                                                                                               \\ \hline\hline

                \multicolumn{13}{l}{Chebyshev-5+hGS-1}                                                                                                                                                                                                                                                                                            \\ \hline
                20                                         & 39(3.6)                                     & 33(3.3)                                     & 33(3.7)                                      & 35(11.1)                       & 35(8.7)     & 29(14.2)  & 34(29.7)   & 36(40.6)    & 35(42.3)   & 34(162.3)  & 28(84.0)    & 34(191.2)   \\
                70                                         & 41(199.9)                                   & 34(173.7)                                   & 33(136.4)                                    & 45(218.5)                      & 35(190.5)   & 35(299.8) & 46(841.0)  & 36(670.3)   & 35(903.4)  & 45(2093.8) & 35(2917.1)  & 34(4099.5)  \\
                84                                         & 41(460.9)                                   & 34(238.3)                                   & 33(325.3)                                    & 45(475.0)                      & 35(529.5)   & 35(610.9) & 44(1332.7) & 36(1028.1)  & 35(1708.7) & 44(4940.4) & 35(4008.7)  & 34(5843.5)  \\ \hline

                \multicolumn{13}{l}{Chebyshev-10+hGS-1}                                                                                                                                                                                                                                                                                           \\ \hline
                20                                         & 34(2.8)                                     & 26(2.2)                                     & 26(3.1)                                      & 38(7.1)                        & 30(7.0)     & 30(10.5)  & 40(25.7)   & 31(32.8)    & 31(37.8)   & 42(108.7)  & 31(86.2)    & 31(152.0)   \\
                70                                         & 36(191.4)                                   & 26(139.4)                                   & 26(114.9)                                    & 40(186.6)                      & 30(234.7)   & 30(264.6) & 42(488.8)  & 31(518.1)   & 31(894.2)  & 44(1746.1) & 32(1518.3)  & 31(3608.8)  \\
                84                                         & 35(411.8)                                   & 26(193.3)                                   & 26(283.1)                                    & 40(436.9)                      & 30(354.9)   & 30(558.3) & 42(1221.2) & 31(970.6)   & 31(1529.2) & 43(4988.0) & 32(3801.3)  & 31(5417.9)  \\ \hline

                \multicolumn{13}{l}{Cholesky+hGS-1}                                                                                                                                                                                                                                                                                               \\ \hline
                20                                         & 33(2.6)                                     & 25(2.2)                                     & 25(3.2)                                      & 37(7.9)                        & 27(6.9)     & 27(11.1)  & 39(28.4)   & 29(34.0)    & 29(34.8)   & 39(129.7)  & 29(103.2)   & 29(150.1)   \\
                70                                         & 35(177.5)                                   & 25(127.6)                                   & 25(177.2)                                    & 39(282.8)                      & 29(216.6)   & 27(320.4) & 41(603.6)  & 29(581.4)   & 29(818.1)  & 41(1928.1) & 29(1582.8)  & 29(3003.8)  \\
                84                                         & 35(401.1)                                   & 25(332.4)                                   & 25(407.1)                                    & 39(657.8)                      & 29(445.5)   & 27(665.0) & 41(1156.5) & 29(869.2)   & 29(1706.1) & 41(4696.6) & 29(2848.4)  & 29(22402.2) \\ \hline
            \end{tabular}}}
    \label{tab: in0.1}
\end{table}

\begin{table}[]
    \caption{Results showing the total number of iterations from low-rank preconditioned GMRES and the total CPU times (in seconds) using preconditioner with $\beta=10^{-4}$ , $\sigma=0.4$, and selected spatial ($n_h$) and stochastic ($n_\xi$) degrees of freedom}
    \centering
    {\setlength{\tabcolsep}{3.5pt}\resizebox{\textwidth}{!}{%
            \begin{tabular}{l|lll|lll|lll|lll}
                \hline
                                                           & \multicolumn{3}{l}{\# iter(t)}              & \multicolumn{3}{l}{\# iter(t)}              & \multicolumn{3}{l}{\# iter(t)}               & \multicolumn{3}{l}{\# iter(t)}                                                                                                             \\
                \hline \hline
                \diagbox{$n_\xi$}{$n_h$}                   &
                \multicolumn{3}{l}{289($h=\frac{1}{2^4}$)} & \multicolumn{3}{l}{1089($h=\frac{1}{2^5}$)} & \multicolumn{3}{l}{4225($h=\frac{1}{2^6}$)} & \multicolumn{3}{l}{16641($h=\frac{1}{2^7}$)}                                                                                                                                              \\ \hline\hline
                $n_{\tau}$                                 & 1                                           & $m_{\xi}$+1                                 & $n_A$                                        & 1                              & $m_{\xi}$+1 & $n_A$     & 1           & $m_{\xi}$+1 & $n_A$      & 1           & $m_{\xi}$+1 & $n_A$      \\ \hline\hline

                \multicolumn{13}{l}{$\sigma_a=0.4$}                                                                                                                                                                                                                                                                                                \\ \hline\hline

                \multicolumn{13}{l}{Chebyshev-5+hGS-1}                                                                                                                                                                                                                                                                                             \\ \hline
                20                                         & 59(5.7)                                     & 36(3.6)                                     & 27(4.5)                                      & 84(13.3)                       & 39(10.1)    & 29(13.2)  & 88(50.9)    & 38(37.8)    & 36(60.4)   & 88(398.2)   & 38(163.3)   & 28(203.0)  \\
                70                                         & 100(275.7)                                  & 38(198.3)                                   & 34(140.1)                                    & 112(498.4)                     & 42(299.3)   & 37(316.8) & 118(1792.6) & 41(697.2)   & 36(929.6)  & 117(7712.6) & 41(2655.2)  & 36(3742.7) \\
                84                                         & 86(720.2)                                   & 36(414.0)                                   & 34(463.0)                                    & 96(986.0)                      & 40(596.4)   & 37(624.8) & 102(3051.9) & 40(1452.6)  & 37(1838.1) & 102(9946.2) & 41(4576.9)  & 36(7289.7) \\ \hline

                \multicolumn{13}{l}{Chebyshev-10+hGS-1}                                                                                                                                                                                                                                                                                            \\ \hline
                20                                         & 69(6.4)                                     & 30(2.8)                                     & 28(3.3)                                      & 77(13.9)                       & 33(8.2)     & 30(9.9)   & 81(50.5)    & 35(34.4)    & 31(38.1)   & 83(208.6)   & 36(148.7)   & 33(162.3)  \\
                70                                         & 92(278.4)                                   & 32(96.5)                                    & 28(125.9)                                    & 105(508.4)                     & 36(175.3)   & 30(269.3) & 110(1318.5) & 37(781.9)   & 32(857.9)  & 110(5212.5) & 38(2920.1)  & 33(4546.2) \\
                84                                         & 79(931.3)                                   & 30(355.0)                                   & 28(411.9)                                    & 89(1373.9)                     & 34(529.2)   & 30(531.0) & 95(3528.4)  & 36(1337.8)  & 32(1576.4) & 97(11001.9) & 38(4413.7)  & 33(5769.8) \\ \hline

                \multicolumn{13}{l}{Cholesky+hGS-1}                                                                                                                                                                                                                                                                                                \\ \hline
                20                                         & 69(6.8)                                     & 29(3.0)                                     & 27(4.0)                                      & 77(17.5)                       & 31(8.2)     & 29(11.7)  & 81(57.5)    & 33(35.2)    & 29(50.3)   & 83(279.8)   & 33(180.9)   & 31(181.8)  \\
                70                                         & 31(159.1)                                   & 31(156.3)                                   & 29(350.5)                                    & 105(730.2)                     & 33(238.6)   & 29(331.6) & 109(1903.1) & 35(726.1)   & 31(953.3)  & 109(7600.0) & 35(2737.8)  & 31(3350.7) \\
                84                                         & 33(753.1)                                   & 33(353.2)                                   & 27(425.5)                                    & 89(1325.8)                     & 33(493.8)   & 29(634.8) & 95(3122.6)  & 35(1055.9)  & 31(1837.7) & 97(10923.5) & 35(4023.8)  & 31(6429.3) \\ \hline
            \end{tabular}}}
    \label{tab: in0.4}
\end{table}

Table~\ref{tab: deter beta} examines the sensitivity to the regularization parameter $\beta$, which balances the tracking term and control cost in the objective functional~\eqref{eq:problem-stoch}.
As $\beta$ decreases from $10^{-2}$ to $10^{-5}$, the optimization problem becomes increasingly dominated by the tracking term.
The iteration counts remain remarkably stable across this range, demonstrating that the hierarchical preconditioner effectively handles varying parameter regimes without requiring problem-specific tuning.
The $n_\tau=m_\xi+1$ truncation consistently performs comparably to the full expansion while maintaining reduced computational cost.
\begin{table}[!ht]
    \caption{Results using the preconditioner with $m_{\xi}$=3 p=3, $\sigma_{\Bbbk}=0.2$, $\beta \in\{10^{-2},10^{-3},10^{-4},10^{-5}\}$ and $n_h=1089(h=\frac{1}{2^5})$.}

    \centering
    {\setlength{\tabcolsep}{3.5pt}\resizebox{\textwidth}{!}{%
            \begin{tabular}{l|lll|lll|lll}
                \hline
                                & \multicolumn{3}{l}{\# iter(t)} & \multicolumn{3}{l}{\# iter(t)} & \multicolumn{3}{l}{\# iter(t)}                                                                             \\
                \hline \hline
                $n_\xi$         & \multicolumn{3}{l}{20}         & \multicolumn{3}{l}{70}         & \multicolumn{3}{l}{84}                                                                                     \\ \hline\hline
                $n_{\tau}$      & 1                              & $m_{\xi}$+1                    & $n_A$                          & 1         & $m_{\xi}$+1 & $n_A$     & 1         & $m_{\xi}$+1 & $n_A$     \\ \hline\hline

                \multicolumn{10}{l}{Chebyshev-5+hGS-1}                                                                                                                                                         \\ \hline
                $\beta=10^{-2}$ & 52(10.0)                       & 32(7.9)                        & 30(13.5)                       & 60(272.3) & 32(149.2)   & 30(363.0) & 56(596.9) & 32(338.9)   & 30(626.9) \\
                $\beta=10^{-3}$ & 52(10.3)                       & 34(8.3)                        & 34(15.5)                       & 60(273.8) & 35(161.7)   & 34(411.4) & 56(594.0) & 34(360.7)   & 34(709.7) \\
                $\beta=10^{-4}$ & 43(11.7)                       & 37(7.6)                        & 35(15.6)                       & 62(279.7) & 37(172.1)   & 35(299.2) & 58(629.0) & 37(393.4)   & 35(731.8) \\
                $\beta=10^{-5}$ & 41(11.2)                       & 38(7.9)                        & 37(17.0)                       & 60(270.6) & 38(175.9)   & 37(317.9) & 56(589.0) & 38(413.3)   & 38(795.5) \\ \hline

                \multicolumn{10}{l}{Chebyshev-10+hGS-1}                                                                                                                                                        \\ \hline
                $\beta=10^{-2}$ & 48(8.7)                        & 30(6.1)                        & 29(11.9)                       & 56(420.8) & 29(147.9)   & 29(359.7) & 52(577.6) & 34(370.8)   & 24(653.1) \\
                $\beta=10^{-3}$ & 48(8.8)                        & 31(6.4)                        & 30(12.6)                       & 56(432.3) & 31(151.6)   & 30(406.9) & 52(554.8) & 31(339.6)   & 30(675.8) \\
                $\beta=10^{-4}$ & 48(8.8)                        & 30(7.1)                        & 30(19.4)                       & 56(419.0) & 31(151.2)   & 30(390.6) & 52(575.0) & 31(339.6)   & 30(676.1) \\
                $\beta=10^{-5}$ & 47(8.7)                        & 30(7.0)                        & 30(13.2)                       & 54(406.4) & 30(146.0)   & 30(383.9) & 51(554.0) & 30(333.0)   & 30(675.1) \\ \hline

                \multicolumn{10}{l}{Chol+hGS-1}                                                                                                                                                                \\ \hline
                $\beta=10^{-2}$ & 47(12.4)                       & 29(6.9)                        & 27(11.9)                       & 55(380.1) & 29(214.4)   & 27(354.4) & 51(975.8) & 29(483.9)   & 27(595.1) \\
                $\beta=10^{-3}$ & 47(12.0)                       & 29(6.4)                        & 27(11.9)                       & 55(396.0) & 29(215.9)   & 27(364.5) & 51(876.9) & 29(455.6)   & 27(596.0) \\
                $\beta=10^{-4}$ & 47(10.4)                       & 29(7.4)                        & 27(12.2)                       & 55(387.9) & 29(216.0)   & 27(343.9) & 53(847.9) & 29(432.8)   & 27(595.5) \\
                $\beta=10^{-5}$ & 45(10.0)                       & 29(6.8)                        & 27(12.0)                       & 53(369.6) & 29(205.5)   & 27(359.3) & 51(829.3) & 29(438.9)   & 27(594.1) \\ \hline
            \end{tabular}}}
    \label{tab: deter beta}
\end{table}

\subsection{Time-dependent case}
\label{subsec:time-num}
This subsection evaluates the all-at-once preconditioner~\eqref{eq: precon for time} for time-dependent optimal control problems.
The discretization results in KKT systems of dimension $n_t \times n_\xi \times n_h$, where $n_t$ denotes the number of time steps.
We investigate the scalability with respect to:
(i)~mesh refinement (Table~\ref{tab: parabolic mesh}),
(ii)~regularization parameter $\beta$ (Table~\ref{tab: parabolic beta}),
(iii)~variance $\sigma_\Bbbk$ (Table~\ref{tab: parabolic CoV}),
(iv)~time-step size discretization $\tau$ (Table~\ref{tab: parabolic Nt}), and
(v)~stochastic dimension $(m_\xi,p)$ (Table~\ref{tab: parabolic KLP}).
Based on the steady-state findings, we employ the Chebyshev-5+hGS-1 configuration unless otherwise noted, reporting results for both stringent ($10^{-6}$) and moderate ($10^{-4}$) tolerances to illustrate practical convergence behavior.
As in the steady-state case, we use $N=1$ in the Richardson iteration (Algorithm~\ref{alg:richardson}) for the (3,3)-block.

From our observations in the steady-state problem, a combination of a 5-step Chebyshev smoother with one step of our hGS method achieves a good balance between the computational cost of matrix operations and GMRES iterations; thus, we typically adopt this combination when testing time-dependent cases as well.

% {\color{blue}We summarize the mean and variance of the state and control variables at the 5th time grid point for mesh size $h=2^{-5}$ with stochastic parameters $m_{\xi}=3$, $p=3$, and $\sigma_\Bbbk=0.2$, and record the residual reduction behavior of the flexible GMRES (fgmres) method under the mean-based approach with $n_A=0$ and the hGS method truncated at $n_A=4$ for the same configuration.}

As indicated in Table~\ref{tab: parabolic mesh}, the 5-step Chebyshev smoother yields consistent iteration counts compared to either the 10-step smoother or direct Cholesky decomposition. As the spatial discretization is refined from $h=\frac{1}{2^3}$ to $h=\frac{1}{2^6}$, representing a growth from 116,640 to 6,084,000 DoF, the iteration count for $n_\tau=m_\xi+1$ exhibits sub-linear growth consistent with the near mesh-independence predicted by the spectral theory in Section~\ref{sec: spectral analysis}.
The $n_\tau=m_\xi+1$ truncation achieves iteration counts comparable to the full expansion ($n_\tau=n_A$) while reducing the cost of assembling and applying the (3,3)-block preconditioner---a trade-off that becomes increasingly favorable as the problem dimension grows.
\begin{table}[!ht]
    \centering
    \caption{Results using hGS method with different truncation settings $n_\tau$ with the time-dependent model for different stopping tolerance and mesh size with $m_{\xi}$=3, $p=3$ $\sigma_\Bbbk$=0.2, $\beta=10^{-4}$ and number of time steps=8 ($\tau=\frac{1}{2^3}$).}
    {
        \setlength{\tabcolsep}{3.5pt}
        \footnotesize
        \begin{tabular}{l|l|lll|lll}\hline %\hline
                            &                         & \multicolumn{3}{l}{tol=$10^{-6}$} & \multicolumn{3}{l}{tol=$10^{-4}$}                                                   \\ \hline
            $h$             & \diagbox{DoF}{$n_\tau$} & 1                                 & $m_{\xi}$+1                       & $n_A$     & 1         & $m_{\xi}$+1 & $n_A$     \\ \hline
            $\frac{1}{2^3}$ & 116,640                 & 57(11.9)                          & 39(8.8)                           & 39(15.7)  & 41(6.8)   & 31(5.7)     & 29(8.8)   \\ \hline
            $\frac{1}{2^4}$ & 416,160                 & 75(39.3)                          & 47(29.6)                          & 45(39.6)  & 55(21.7)  & 35(14.9)    & 35(24.7)  \\ \hline
            $\frac{1}{2^5}$ & 1,568,160               & 83(150.9)                         & 53(103.3)                         & 51(158.6) & 66(90.6)  & 43(58.8)    & 39(92.1)  \\ \hline
            $\frac{1}{2^6}$ & 6,084,000               & 84(570.8)                         & 55(383.3)                         & 53(585.1) & 68(401.4) & 42(252.2)   & 42(424.1) \\ \hline
            % & 
        \end{tabular}}
    \label{tab: parabolic mesh}
\end{table}

Table~\ref{tab: parabolic beta} examines four orders of magnitude for $\beta$, ranging from $10^{-2}$ (control-dominant) to $10^{-8}$ (tracking-dominant).
The mean-based preconditioner ($n_\tau=1$) exhibits strong dependence on $\beta$, with iteration counts decreasing as $\beta$ decreases (since smaller $\beta$ yields problems dominated by the simpler tracking term).
In contrast, the $n_\tau=m_\xi+1$ truncation maintains stable iteration counts across all tested values, demonstrating that the hGS preconditioner automatically adapts to the problem structure without manual parameter tuning.
This robustness confirms the applicability of the theoretical framework across diverse parameter regimes.
\begin{table}[!ht]
    \centering
    \caption{Results using hGS method with different truncations setting $n_\tau$ with the time-dependent model for different stopping tolerance and $\beta$, mesh size $h=\frac{1}{2^5}$, $m_{\xi}$=3, $p=3$, $\sigma_\Bbbk$=0.2, and number of time steps=8 ($\tau=\frac{1}{2^3}$), which results in 1,568,160 DoF.}
    {
        \setlength{\tabcolsep}{3.5pt}
        \footnotesize
        \begin{tabular}{l|lll|lll}\hline
                                        & \multicolumn{3}{l}{tol=$10^{-6}$} & \multicolumn{3}{l}{tol=$10^{-4}$}                                                   \\ \hline
            \diagbox{$\beta$}{$n_\tau$} & 1                                 & $m_{\xi}$+1                       & $n_A$     & 1         & $m_{\xi}$+1 & $n_A$     \\ \hline
            $10^{-2}$                   & 107(206.1)                        & 69(136.3)                         & 69(214.8) & 74(102.0) & 48(70.4)    & 47(112.6) \\ \hline
            $10^{-3}$                   & 116(164.4)                        & 74(109.4)                         & 74(174.6) & 66(90.6)  & 42(60.7)    & 42(100.0) \\ \hline
            $10^{-6}$                   & 97(137.0)                         & 77(115.6)                         & 77(182.5) & 57(77.4)  & 45(65.3)    & 44(103.7) \\ \hline
            % $10^{-8}$&78(109.2)&76(114.7)&77(183.1)&44(60.3)&42(61.1)&42(99.2)\\ \hline
        \end{tabular}
    }
    \label{tab: parabolic beta}
\end{table}

Table~\ref{tab: parabolic CoV} explores the range $\sigma_\Bbbk \in \{0.01, 0.02, 0.05, 0.1, 0.2, 0.4\}$, spanning from nearly deterministic to highly uncertain regimes.
The mean-based preconditioner ($n_\tau=1$) exhibits significant degradation as uncertainty increases, whereas the~$n_\tau=m_\xi+1$ truncation maintains stable iteration counts across the entire range. Especially when $\sigma_\Bbbk$ increases from 0.2 to 0.4, the mean-based preconditioner performs poorly with a large number of iterations, but the hGS method maintains robust performance.
This robustness confirms that the hierarchical preconditioner effectively captures the essential stochastic structure without requiring full expansion of all coupling coefficients.
\begin{table}[!ht]
    \centering
    \caption{Results using hGS method with different truncations setting $n_\tau$ with the time-dependent model for different stopping tolerance and $\sigma_\Bbbk$, $\beta=10^{-4}$, mesh size $n_h=1089$, $m_{\xi}$=3, $p=3$, and number of time steps=8 ($\tau=\frac{1}{2^3}$), which results in 1,568,160 DoF.}
    {\setlength{\tabcolsep}{3.5pt}\footnotesize
        \begin{tabular}{l|lll|lll}\hline
                                               & \multicolumn{3}{l}{tol=$10^{-6}$} & \multicolumn{3}{l}{tol=$10^{-4}$}                                                    \\ \hline
            \diagbox{$\sigma_\Bbbk$}{$n_\tau$} & 1                                 & $m_{\xi}$+1                       & $n_A$     & 1          & $m_{\xi}$+1 & $n_A$     \\ \hline
            0.01                               & 51(77.7)                          & 51(85.3)                          & 51(133.9) & 43(68.6)   & 41(72.7)    & 41(112.4) \\ \hline
            0.02                               & 53(80.8)                          & 51(82.2)                          & 51(130.0) & 43(70.2)   & 41(71.1)    & 41(113.7) \\ \hline
            0.05                               & 59(90.2)                          & 51(82.1)                          & 51(130.9) & 47(78.3)   & 41(72.6)    & 41(113.1) \\ \hline
            0.1                                & 65(99.2)                          & 51(82.4)                          & 51(130.4) & 53(86.2)   & 42(70.6)    & 41(116.5) \\ \hline
            0.2                                & 83(150.9)                         & 53(103.3)                         & 51(158.6) & 66(90.6)   & 43(58.8)    & 39(92.1)  \\ \hline%done
            0.4                                & 129(200.9)                        & 57(92.3)                          & 51(130.5) & 100(179.3) & 47(82.3)    & 43(119.2) \\ \hline
        \end{tabular}}
    \label{tab: parabolic CoV}
\end{table}
Table~\ref{tab: parabolic Nt} investigates the all-at-once system scalability by varying $n_t$ from 4 to 256 (time steps $\tau \in \{1/4, 1/16, 1/64, 1/256\}$), corresponding to total system sizes ranging from 784,080 to over 12.5 million DoF.
As the temporal resolution increases, the coupled space-time-stochastic system grows proportionally, yet the $n_\tau=m_\xi+1$ truncation maintains sub-linear iteration growth relative to system size.
The computational time scales approximately linearly with DoF, confirming the efficiency of the all-at-once preconditioner for massively coupled systems.
\begin{table}[!ht]
    \centering
    \caption{Results using hGS method with different truncations setting $n_\tau$ with the time-dependent model for different stopping tolerance and $\sigma_\Bbbk$, $\beta=10^{-4}$, mesh size $n_h=1089$, $m_{\xi}=3,\, p=3$, and number of time steps=8 ($\tau=\frac{1}{2^3}$).}
    {\setlength{\tabcolsep}{3.5pt}\footnotesize
        \begin{tabular}{l|l|lll|lll}\hline
                    &                         & \multicolumn{3}{l}{tol=$10^{-6}$} & \multicolumn{3}{l}{tol=$10^{-4}$}                                                     \\ \hline
            $\tau$  & \diagbox{DoF}{$n_\tau$} & 1                                 & $m_{\xi}$+1                       & $n_A$      & 1         & $m_{\xi}$+1 & $n_A$      \\ \hline
            $1/2^2$ & 784,080                 & 81(78.6)                          & 51(51.1)                          & 51(71.1)   & 66(53.9)  & 43(37.6)    & 42(57.2)   \\ \hline
            $1/2^4$ & 3,136,320               & 85(238.2)                         & 55(166.4)                         & 53(258.7)  & 68(174.2) & 44(124.3)   & 43(203.7)  \\ \hline
            $1/2^6$ & 12,545,280              & 101(1029.8)                       & 67(733.3)                         & 65(1205.8) & 78(840.5) & 53(604.2)   & 52(1028.6) \\ \hline
            % $1/2^8$&5575680($256\times 20\times 1089$)&181(7391.9)&&129(9584.8)&138(5717.4)&108(4862.3)&108(8225.1)\\ \hline
            % $1/2^{10}$&22302720&103(141.8)&63(91.5)&63(150.8)&66(90.6)&43(58.8)&39(92.1)\\ \hline
        \end{tabular}}
    \label{tab: parabolic Nt}
\end{table}

Finally, Table~\ref{tab: parabolic KLP} compares three gPC configurations:~$(m_\xi, p) \in \{(3,3), (4,4), (6,3)\}$, yielding $n_\xi \in \{20, 70, 84\}$ basis functions.
Table~\ref{tab: parabolic KLP} varies the stochastic discretization parameters $(m_\xi, p)$, exploring both the number of random variables and polynomial order.
As $n_\xi$ increases from 20 to 84, the iteration count for $n_\tau=m_\xi+1$ grows modestly, demonstrating near-independence from the stochastic discretization level.
This behavior confirms the effectiveness of the hierarchical truncation strategy in maintaining spectral properties across varying gPC expansion settings.
\begin{table}[!ht]
    \centering
    \caption{Results using hGS method with different truncations setting $n_\tau$ with the time-dependent model for different stopping tolerance and stochastic setting, $\beta=10^{-4}$, mesh size $n_h=1089$, $\sigma_\Bbbk$=0.2, and number of time steps=8 ($\tau=\frac{1}{2^3}$).}
    {\setlength{\tabcolsep}{3.5pt}\footnotesize
        \begin{tabular}{l|l|lll|lll}\hline
                    &                         & \multicolumn{3}{l}{tol=$10^{-6}$} & \multicolumn{3}{l}{tol=$10^{-4}$}                                                      \\ \hline
            $n_\xi$ & \diagbox{DoF}{$n_\tau$} & 1                                 & $m_{\xi}$+1                       & $n_A$      & 1          & $m_{\xi}$+1 & $n_A$      \\ \hline
            20      & 1,568,160               & 83(150.9)                         & 53(103.3)                         & 51(158.6)  & 66(90.6)   & 43(58.8)    & 39(92.1)   \\ \hline %done
            70      & 5,488,560               & 95(2435.3)                        & 53(1397.4)                        & 52(2964.7) & 74(1809.2) & 43(1044.9)  & 43(2476.2) \\ \hline
            84      & 6,586,272               & 93(4092.4)                        & 56(2670.3)                        & 52(5207.4) & 70(3059.4) & 43(2107.7)  & 43(4413.4) \\ \hline
        \end{tabular}}
    \label{tab: parabolic KLP}
\end{table}

The time-dependent experiments establish that the proposed all-at-once preconditioner maintains robust performance across a wide range of problem parameters and discretization levels.
Four key findings emerge from these results.
First, regarding \textit{near mesh-independence}, iteration growth remains sub-linear with spatial refinement (Table~\ref{tab: parabolic mesh}), consistent with the spectral bounds derived in Section~\ref{sec: spectral analysis}.
Second, concerning \textit{parameter robustness}, the hGS method adapts automatically to varying $\beta$ (Table~\ref{tab: parabolic beta}) and $\sigma_\Bbbk$ (Table~\ref{tab: parabolic CoV}) without manual tuning.
Third, regarding \textit{truncation efficiency}, the $n_\tau=m_\xi+1$ strategy consistently delivers performance comparable to the full expansion at significantly reduced cost. These results demonstrate that the hierarchical preconditioning framework extends seamlessly from steady-state to time-dependent problems, providing a practical and theoretically-grounded solution for large-scale stochastic optimal control problems.

\section{Conclusions}
\label{sec:conclusion}
In the paper, we designed, analyzed, and implemented a novel hierarchical preconditioning strategy for large-scale stochastic optimal control problems. Our approach leverages a truncated stochastic expansion within a block-structured preconditioner for the Karush-Kuhn-Tucker (KKT) system, striking an effective balance between computational cost and preconditioning quality. Numerical results confirm that the proposed hGS method consistently outperforms both standard mean-based preconditioner and computationally intensive full-expansion preconditioner across a wide range of problem parameters.

A key contribution of this work is the extension of the truncated gPC framework to time-dependent problems. We developed and tested a tailored hGS preconditioner within an all-at-once discretization scheme, demonstrating the versatility and effectiveness of our approach for these challenging, large-scale scenarios. Comprehensive numerical experiments on benchmark problems have validated the robustness and numerical efficiency of the proposed algorithms. %Future research could involve extending this preconditioning framework to problems with more complex PDE constraints, such as the Navier-Stokes equations, or investigating its application to optimal control problems with inequality constraints.
\bibliographystyle{siamplain}
%\bibliography{control_prec.bib}
\bibliography{ref.bib}
\clearpage
\appendix
\section{Auxiliary Results}
\label{sec:appendix}

\subsection{Proof of Lemma~\ref{lem:error_bound}}
\begin{proof}
    Consider the Rayleigh quotient:
    \[
        \frac{\mathbf{v}^T A \mathbf{v}}{\mathbf{v}^T B \mathbf{v}} = \frac{\mathbf{v}^T (B + E) \mathbf{v}}{\mathbf{v}^T B \mathbf{v}} = 1 + \frac{\mathbf{v}^T E \mathbf{v}}{\mathbf{v}^T B \mathbf{v}}.
    \]
    Applying the assumption $|\frac{\mathbf{v}^T E \mathbf{v}}{\mathbf{v}^T B \mathbf{v}}| \le \delta$ directly yields the bounds.
\end{proof}

\subsection{Proof of Lemma~\ref{lem:congruence}}
\begin{proof}
    Let $(\lambda, \mathbf{x})$ be an eigenpair satisfying the generalized eigenvalue problem $C\mathbf{x} = \lambda D\mathbf{x}$, with eigenvector $\mathbf{x} \neq \mathbf{0}$. We perform a change of variables by setting $\mathbf{x} = Q\mathbf{y}$. Since $Q$ is nonsingular, $\mathbf{x} \neq \mathbf{0}$ implies that the transformed vector $\mathbf{y} \neq \mathbf{0}$.
    Substituting $\mathbf{x} = Q\mathbf{y}$ into the original problem gives:
    \[
        C(Q\mathbf{y}) = \lambda D(Q\mathbf{y}).
    \]
    Multiplying from the left by $Q^T$, we obtain:
    \[
        (Q^T C Q)\mathbf{y} = \lambda (Q^T D Q)\mathbf{y}.
    \]
    This final expression is the generalized eigenvalue problem for the pair~$(Q^T C Q, Q^T D Q)$, which is satisfied by the same eigenvalue $\lambda$ with the transformed eigenvector $\mathbf{y}$. Therefore, the sets of eigenvalues for both pairs are identical.
\end{proof}
\subsection{Proof of Lemma~\ref{lem:sig<=l2}}
\begin{proof}
    \zhl{
        Let the SVD of $A$ be $A = U\Sigma V^*$, where $U, V$ are unitary and $\Sigma = \operatorname{diag}(\sigma_1, \ldots, \sigma_n)$ with $\sigma_1 \geq \cdots \geq \sigma_n = \sigma_{\min}(A) \geq 0$.

        For any nonzero vector $\mathbf{v}$, let $\mathbf{w} = V^*\mathbf{v}$. Since $V$ is unitary, $\|\mathbf{w}\|_2 = \|\mathbf{v}\|_2$. Then
        \[
            \|A\mathbf{v}\|_2^2
            = \|U\Sigma V^*\mathbf{v}\|_2^2
            = \|\Sigma \mathbf{w}\|_2^2
            = \sum_{i=1}^{n} \sigma_i^2 |w_i|^2
            \geq \sigma_{\min}^2 \sum_{i=1}^{n} |w_i|^2
            = \sigma_{\min}^2 \|\mathbf{w}\|_2^2
            = \sigma_{\min}^2 \|\mathbf{v}\|_2^2.
        \]
        Taking square roots on both sides gives $\|A\mathbf{v}\|_2 \geq \sigma_{\min}(A)\|\mathbf{v}\|_2$.}
\end{proof}

\subsection{Proof of Lemma~\ref{lemma ineq}}
\begin{proof}
    From the triangle inequality, for any vector $x$, with $\|x\|_2=1$, we have
    \[
        \|(A+B)x\|_2 = \|Bx - (-A)x\|_2 \geq \|Bx\|_2 - \|Ax\|_2.
    \]
    Taking the minimum over all unit vectors $x$ on both sides of the inequality, we get
    \[
        \min_{\|x\|_2=1} \|(A+B)x\|_2 \geq \min_{\|x\|_2=1} \left( \|Bx\|_2 - \|Ax\|_2 \right).
    \]
    Using the property that $\min(f-g) \geq \min(f) - \max(g)$, we obtain
    \[
        \min_{\|x\|_2=1} \left( \|Bx\|_2 - \|Ax\|_2 \right) \geq \min_{\|x\|_2=1} \|Bx\|_2 - \max_{\|x\|_2=1} \|Ax\|_2.
    \]
    By the definitions of the minimum singular value and the operator norm, the above expression is equivalent to $\sigma_{\min}(A+B)\geq \sigma_{\min}(B)-\|A\|_2$.
\end{proof}

\end{document}